\documentclass[twoside,12pt]{article}
\usepackage{amsmath,amscd,amsthm,amssymb,amsxtra,latexsym,epsfig,epic,graphics}
\usepackage[matrix,arrow,curve]{xy}
\usepackage{graphicx}
\usepackage{diagrams}

\def\antiddot{\mathinner{\mkern1mu\raise1pt\vbox{\kern7pt\hbox{.}}\mkern2mu
        \raise4pt\hbox{.}\mkern2mu\raise7pt\hbox{.}\mkern1mu}}

\newcommand{\AAA}{{\mathbb A}}

\newcommand{\FF}{{\mathbb F}}

\newcommand{\PP}{{\mathbb P}}

\newcommand{\RR}{{\mathbb R}}

\newcommand{\ZZ}{{\mathbb Z}}
\newcommand{\LF}{{\rm\bf L}}
\newcommand{\RF}{{\rm\bf R}}
\newcommand{\TF}{{\rm\bf T}}
\newcommand{\UF}{{\rm\bf U}}

\newcommand{\HH}{{\rm{H}}}

\newcommand{\Ext}{{\rm{Ext}}}

\newcommand{\rad}{{\rm{rad\;}}}

\newcommand{\id}{{\rm{id}}}

\newcommand{\coker}{{\rm{coker}\,}}

\newcommand{\s}{\mathcal}

\newcommand{\sB}{{\s B}}

\newcommand{\sE}{{\s E}}
\newcommand{\sF}{{\s F}}
\newcommand{\sG}{{\s G}}
\newcommand{\sH}{{\s H}}

\newcommand{\sK}{{\s K}}

\newcommand{\sO}{{\s O}}

\newcommand{\tensor}{\otimes}

\newcommand{\punkt}{\hspace{-.3ex}\raise.15ex\hbox to1ex{\Huge.}}

\newcommand{\fix}{{\bf **** fix: }}

\DeclareMathOperator{\Sym}{Sym}
\DeclareMathOperator{\reg}{reg}
\DeclareMathOperator{\Spec}{Spec}
\DeclareMathOperator{\Proj}{Proj}
\DeclareMathOperator{\Hom}{Hom}
\DeclareMathOperator{\sHom}{\sH om}

\DeclareMathOperator{\pd}{pd}

\DeclareMathOperator{\Tor}{Tor}

\DeclareMathOperator{\rank}{rank}

\newcommand{\gm}{\mathfrak m}

\newcommand{\Mac}{{\texttt {MACAULAY2}}}

\newtheorem{theorem}{Theorem}[section]
\newtheorem{lemma}[theorem]{Lemma}
\newtheorem{proposition}[theorem]{Proposition}
\newtheorem{corollary}[theorem]{Corollary}
\newtheorem{conjecture}{Conjecture}[section]
\theoremstyle{definition}

\newtheorem{remark}[theorem]{Remark}

\newtheorem{example}[theorem]{Example}

\date{September 27, 2006}
\title{Relative Beilinson Monad and Direct Image for Families
of Coherent Sheaves} 
\author{David Eisenbud and Frank-Olaf Schreyer}

\begin{document}

\maketitle

\begin{abstract}

The higher direct image complex of a coherent sheaf (or finite complex
of coherent sheaves) under a projective morphism is a fundamental
construction that can be defined via a \v Cech complex or an injective
resolution, both inherently infinite constructions. Using
free resolutions it can be defined in finite terms. Using exterior
algebras and relative versions of theorems of Beilinson and
Bernstein-Gel'fand-Gel'fand, we give an alternate
and generally more efficient description in
finite terms.

Using this exterior algebra description we can characterize 
the generic finite free complex of a given shape
as the direct image of an easily-described
vector bundle. We can also
give explicit descriptions of the loci
in the base spaces of flat families of sheaves in which some
cohomological conditions are satisfied---for example, the loci where
vector bundles on projective space split in a certain way, or the loci where 
a projective morphism has higher dimensional fibers. 

Our approach is so explicit that it yields an algorithm suited for
computer algebra systems.

\end{abstract}

\section*{Introduction}

 Let $A$ be a Noetherian ring,  let
$\sF$ be a coherent sheaf on $\PP^n_A= \PP^n \times\Spec A$, and let
$\pi$ be the projection map  $\PP^n_A\to \Spec A$:
$$
\xymatrix{ 
\sF & \PP^n_A=\PP^n\times{\Spec A} \ar[d]^\pi \cr
R\pi_* \sF & \Spec A \cr }
$$

A locally free complex $R\pi_* \sF$ of coherent sheaves on ${\Spec A}$ (or, equivalently, of $A$-modules)
is defined up to quasi-isomorphism by pushing forward an injective resolution
of $\sF$ and choosing a bounded complex of finitely generated $A$-modules 
quasi-isomorphic to it. There is also
a description in finite terms using free resolutions over
the polynomial ring $A[x_0,\dots,x_n]$ (see Proposition \ref{elementary push-forward}).
 In this paper we will give a more efficient construction,
using (finite parts of) resolutions over an exterior algebra.
This makes it possible to compute, for example,
the loci in flat families of sheaves where certain cohomological conditions are
satisfied, such as the loci where certain decompositions occur in a family of
vector bundles, or to detect higher dimensional fibers.

For example, if $\sF$ is flat over $A$, and $A$
is local with maximal ideal $\gm$, then by a
Theorem of Grothendieck [EGA,7.7], see also
[Mumford,II.5] or [Hartshorne,III.12],
the complex $R\pi_* \sF \in D^b(A)$ is represented by a
minimal complex
$$0\to A^{h^0} \to A^{h^1} \to \ldots \to A^{h^n} \to 0$$
of free $A$-modules, which is unique up to isomorphism.
Because $\sF$ is flat the formation of this complex commutes with 
base change, and hence
$$h^i=\dim_{K} \HH^i(\sF\tensor  K),$$
where $K=A/\gm$ denotes the residue class field of $A$. 
If $\sF$ is given explicitly, then our techniques
algorithmically compute the matrices in this complex.

To formulate the main result we introduce some notation.
Let $W=\pi_* \sO_{\PP^n_A}(1)$ be the free 
$A$-module of rank $n+1$ underlying $\PP^n_A$, and
let $x_0,\ldots,x_n$ be a basis of $W$. 
The scheme
$\PP^n_A$ is $\Proj S$ where $S=\Sym W \cong A[x_0,\ldots,x_n]$.
Let $M=\sum_d M_d$ be a graded $S$-module whose sheafification is $\sF$,
and let $\reg(M)$ denote its Castelnuovo-Mumford regularity
(as defined, in this relative case, below.)

Let $V=\Hom_A(W,A)$ be the dual of $W$,
 with basis $e_0,\ldots,e_n$ dual to $x_0,\dots,x_n$. 
Let  
$$
E =\Lambda V =\oplus_{i=0}^{n+1} \Lambda^i V
$$
be the exterior algebra on $V$. We give the 
variables $e_i$ degree $-1$. We write $(V)$ for the 
\emph{augmentation ideal} generated by $V$ in $E$.

Any projective $E$-module $M$ can be written in the form $E\otimes_AN$,
where $N$ is the projective $A$-module $M/(V)M$. More 
generally, if $N$ is any $A$-module, we say that $E\otimes N$
is a \emph{relatively projective} $E$-module. We make a 
similar definition for modules over $S$. There is a well-developed
theory of relative homological algebra, including relative
projective resolutions, due to Hochschild [Hoch56].
We review some elementary
facts about it in 
Section \ref{relative projectives}. For us relative projectives will
be useful because modules generally have  smaller
resolutions in terms of relative projectives than in 
terms of projectives.

Corresponding to the graded $S$-module $M$ there is
a complex of $E$-modules
$$
\cdots \to E\otimes M_d \to E\otimes M_{d+1}\to E\otimes M_{d+2}\to\cdots
$$
where $\otimes$ denotes $\otimes_A$,
the module $E\otimes M_d$ is in cohomological degree $d$
and is generated by $M_d$ regarded as an $A$-module
concentrated in degree $d$. The differentials in the complex
are given by
$$
a\otimes m\mapsto \sum_i e_ia\tensor x_im.
$$

\def\quis{{\sim}}
In what follows we will write $\quis$ for quasi-isomorphism
of (bounded above) complexes---the equivalence relation generated by declaring
two complexes quasi-isomorphic if there is a map between them that
is an isomorphism on homology (for free bounded-above complexes
this is the same as homotopy equivalence.)

\begin{theorem}[Main Theorem]\label{main}
Suppose $s\ge \max(0,\reg(M))$, and set
$P^s=\ker (E\otimes M_{s+1} \to E\otimes M_{s+2})$.
If $\TF$ is a graded relatively projective resolution of $P^s$, regarded as
an $E$-module concentrated in cohomological degree $s$, then 
$$
R\pi_* \sF \quis (\TF \tensor_E A)_0.
$$
\end{theorem}

\begin {corollary}
With hypotheses and notation as in Theorem \ref{main},
$$
R^i\pi_* \sF \cong \Tor^E_{s-i}(P^s,A)_0,
$$
\end{corollary}
\begin{proof} Take $\TF$ to be a projective resolution.
\end{proof}

It is interesting to compare the result of Theorem \ref{main} with
the following more elementary construction of 
$
R\pi_* \sF.
$
Here we use a different notion of regularity, derived from
a minimal relatively projective $A$-split resolution (see Section \ref{relative projectives}.
It seems reasonable to conjecture that this bound
can be improved along the lines of Smith [Sm00], who used the
same idea to compute global Ext, and thus the individual
$R^i\pi_*$ functors. However we cannot simply take
 $d>\reg M$ since this would allow us to take $I=S$
 whenever $\reg M \leq n$, leading to a false result.

\begin{theorem}
\label{elementary push-forward} Assume for simplicity $n\geq 1$.
Let $M$ be a finitely generated $A[x_0,\dots,x_n]$ module,
and let $\sF$ be the corresponding sheaf. Let $d$ be such
that $M$ admits a relatively projective $A$-split resolution by
modules $S\otimes N_i$ where each $N_i$ is a graded $A$-module
that is zero in degrees $>d$ (see Theorem \ref{relative syzygy}).
Suppose also that $I\subset S$ is an ideal such that $S/I$ is projective
as an $A$-module, $I$ is contained in $(x_0,\dots,x_n)^{d-n}$,
and $I$ contains some power of each $x_i$. If
$\FF$ is a free resolution of $I$, then
$$
R\pi_* \sF\quis (\Hom_S(\FF,M)_0.
$$
\end{theorem}

Among the examples of ideals $I$ such that $S/I$ is $A$-projective
are all ideals generated by monomials; and any ideal generated by
a homogeneous complete intersection of length $n+1$ that contains
a power of each of the $x_i$.

\noindent We give a proof in section \ref{relative projectives}.
\fix{Put in examples of these methods
---including one were the second where this is very inefficient.}

\bigbreak

In the first section we recall some facts about
relative projectives, and use them to prove
Theorem \ref{elementary push-forward}.

The next three sections of the paper
contain a description of $\TF$ and a proof of 
Theorem \ref{main}.
 Its main ingredients are
\begin{enumerate}
\item An effective construction of a relative Beilinson monad.
\item The fact that the relative Beilinson monad for 
$\sF$ is $\pi_*$-acyclic---that is, $R^i\pi_*$ 
vanishes on all terms of the monads for all $i>0$.
\end{enumerate}
 
In the remainder of the paper we carry out this construction
of the direct image complex in 
three examples. The first concerns the versal deformation
of a rank $r$ vector bundle of the form $\sO^{r-1}\oplus \sO(d)$ on $\PP^1$ 
(this example can, of course be treated
by other means, and we sketch a
more elementary alternative.) The result gives determinantal 
equations for the loci, in the base
space of the deformation, 
of bundles of a given splitting type; we conjecture that these determinantal
equations actually generate the prime ideals corresponding to the
loci in question. The case $r=2$ has a considerable history, and can
already be found,
in equivalent form and
without proof, in Room [1938]. 

Our second example treats the direct images of certain
sheaves on the resolution of an elliptic
singularity; this example seems amenable only to computation.

The last example was the most surprising to us. It leads to a new
description of the variety of complexes, in 
Theorem \ref{all occur 2}, and answers the question, ``Which
complexes appear as direct images of vector bundles?"

\begin{theorem} Let $A$ be a Noetherian ring, and let 
$$ \FF:\quad 0 \to A^{\beta_0} \to  A^{\beta_1} \to \ldots \to A^{\beta_n} \to 0$$
be any finite complex of finitely generated free $A$-modules of length $n$.
There is a vector bundle $\sF$ on $\PP^n_A$ such that $\FF$ represents
$R\pi_*\sF$. 
\end{theorem}
 
In fact the bundle $\sF$ can be given quite explicitly; see the proof of
Theorem \ref{all occur}.

An implementation of the resulting algorithms in the system \Mac\, and all
the
examples treated as illustrations of the main result can be
downloaded from {\tt
http//www.math.uni-sb.de/\~{}ag-schreyer/computeralgebra}.

Acknowledgment: We thank Mike Stillman and Dan Grayson for their
program \Mac\ and for their quick responses to our needs.
The second author thanks the MSRI for financial support and hospitality.

\medbreak
\noindent{\bf Notation:}
Throughout this paper, $A$ will denote a Noetherian ring
and $S$ will denote the polynomial ring $A[x_0,\dots,x_n] =\Sym(W)$,
where $W$ is the free $A$-module of rank $n+1$ generated by the 
$x_i$. We grade $S$ with $A$ in degree 0 and the $x_i$ in degree 1.
There is a canonical projection $\pi: \PP^n_A\to \Spec A$ corresponding
to the inclusion $A\subset S$. We denote by $M$ a
finitely generated graded $S$-module, and by
$\sF= \widetilde M$ the associated sheaf on $\PP^n_A = \Proj S$.
We write
$R\pi_*\sF$ for the (total) direct image, which is a complex
of $A$ modules determined up to quasi-isomorphism. 

We write $V$ for the $A$-dual of $W$; the module $V$ is
thus a free $A$-module concentrated in degree $-1$. We let
$e_0,\dots,e_n$ denote a basis of $V$ dual
to $x_0,\dots x_n$. Let $E$ be the exterior algebra 
$E=\bigwedge V=A\langle e_0,\dots,e_n\rangle$.
If $A$ is a local ring with maximal ideal $\gm$ then $E$ is a local ring
with maximal ideal $\gm_E=\gm E +(e_0,\ldots,e_n)$. 

\section{Relative Projectives; Proof of Theorem \ref{elementary push-forward}}
\label{relative projectives}

The definitions and results on relative projectives 
in this section come from
Hochschild [Hoch56]; we review them
for the convenience of the reader.

Let $A$ be a commutative Noetherian ring, and let $B$ be an
$A$-algebra.
A $B$-module $P$ is said to be \emph{relatively projective} (with
respect to $A$) if for every surjective $M\to M'$ of $B$-modules
which splits as a map of $A$ modules, the induced map
$\Hom_B(P,M)\to \Hom_B(P,M')$
is also surjective; that is, $B$-linear maps from $P$ can be lifted along
$A$-split surjections. Modules of the form
$B\otimes_AN$, where $N$ is any $A$-module,
are called \emph{relatively free}. 

\begin{proposition} \label{elementary}
A $B$-module $P$ is relatively
projective if and only if it is a direct summand of
a relatively free $B$-module.
If $B$ is Noetherian and positively (or negatively) graded,with
degree 0 component $A$, then
every finitely generated 
graded relatively projective $B$-module 
is relatively free.
\end{proposition}

\begin{proof}
If $P=B\otimes_A N$, then $\Hom_B(P,M)=\Hom_A(N,M)$
and the relative projectivity of $P$ follows from the
definition. It is also immediate that that a
direct summand of a relatively projective module is relatively
projective. If $P$ is any $B$-module, then the 
natural surjection $B\otimes_AP\to P$ is $A$-split by the
map sending $p\in P$ to $1\otimes p\in B\otimes_AP$.
Thus if $P$ is relatively projective it is a summand of
the relatively free module $B\otimes_AP$.

Now suppose that $B$ is graded, and 
$P$ is a graded relatively projective $B$-module.
Since the map $B\to B/B_+=A$ is $A$-split,
 we may write $P/B_+P$ as an $A$-submodule of $P$.
The induced
map of $B$-modules $B\otimes_A (P/B_+P) \to P$ is
surjective by Nakayama's Lemma, and is $A$-split. 
Since $P$ is relatively projective, we
can lift the identity map of $P$ and get a splitting
$P\to B\otimes_A (P/B_+P)$. The complementary
summand goes to zero under the map
$B\otimes_A (P/B_+P) \to P$. But this map becomes an isomorphism modulo
$B_+$. Using Nakayama's Lemma again, we see that the kernel
is zero.
\end{proof}

We now return to our basic case, where $S$ and $E$ are the
exterior and symmetric algebras over $A$ as in the 
introduction.

If $T$ is a graded relatively projective $E$-module,
then by Proposition \ref{elementary} we can write
$T=E\otimes_A N$, 
where $N=\bigoplus_j N_j$ is a graded $A$-module,
with $N_j$ concentrated in degree $j$.
Here we regard $A$ as concentrated in degree 0. We think of $N_j$
as the module of ``generators of $T$ of degree $j$".
A \emph{relatively projective resolution} of a module $M$ is
defined to be a resolution by relatively projective modules.
This terminology differs a little from Hochschild's: he adds
the requirement that the resolution be split exact as a sequence
of $A$-modules. We will encounter both sorts of resolutions below,
so it seems simpler to require $A$-split exactness explicitly
when we need it.

Relative projective resolutions are often  smaller than
projective resolutions. For example, if $A$ is not a regular
ring, then some finitely generated modules over the polynomial ring
$S=A[x_0,\dots,x_n]$ do not have finite projective resolutions, but
they do have finite relatively projective resolutions:

\begin{theorem}\label{relative syzygy}
[Hoch56, Theorem 3]
Every finitely generated 
(graded) $S$-module $M$ has a finite (graded) $A$-split
relatively
projective resolution, of length at most $n+1$. \qed
\end{theorem}

We will use Theorem \ref{relative syzygy} together
with the following Lemma to prove
Theorem \ref{elementary push-forward}.

\begin{lemma}\label{hom-tensor}
If $X$ is an $S$-module and $N$ is an $A$-module
there is a natural transformation
$\Hom_S(X,S)\otimes_A N \to \Hom_S(X, S\otimes_A N)$
that is an isomorphism whenever $X$ is $A$-projective.
\end{lemma}

\begin{proof}

The map is the composition
$
\Hom_S(X,S)\otimes_A N\cong
\Hom_S(X,S)\otimes_S (S\otimes_AN) \cong
\Hom_S(X,S)\otimes_S \Hom_S(S, S\otimes_AN) \to
\Hom_S(X, S\otimes_AN),
$
where the last map is given by composition of homomorphisms.
If $N$ is a free $A$-module, the map is trivially an isomorphism.
Thus
if  $\phi: G_1\to G_0$ is an $S$-free presentation of $X$, then
we get an isomorphism
$\Hom_S(\phi,S)\otimes_A N\to \Hom_S(\phi,S\otimes_AN)$.
The kernel of the map $\Hom_S(\phi,S\otimes_AN)$ is
$
\Hom_S(X, S\otimes_AN).
$
 If $X$ is $A$-projective then $\phi$ is
$A$-split, so the kernel of 
$\Hom_S(\phi,S)\otimes_A N$
is equal to 
$(\ker \Hom_S(\phi,S))\otimes_A N = \Hom_S(X,S)\otimes_A N$
as required.
\end{proof}

\noindent{\it Proof of Theorem \ref{elementary push-forward}}
First, suppose that $M=S(-d)$, the free graded module of rank 1
generated in degree $d$, which has regularity $d$. In this case 
$\sF=\sO_{\PP^n_A}(-d)$.
At most one $R^i\pi_*\sO_{\PP^n_A}(-d)$ is nonzero, 
so $R\pi_*\sF$ is quasi-isomorphic
to this nonzero module, as a complex concentrated in cohomological degree $i$.

Our hypotheses on $I$ imply that the projective dimension of $I$ is
$n$, and that
$$
H^i\Hom_S(\FF, M)=
\begin{cases}
\Hom_S(I, S(-d)) = S(-d) & \hbox{if $i=0$;} \\
\Hom_A(S/I, A)(-d+n+1) &\hbox{if $i=n$;}\\
0 & \hbox{otherwise.}\\
\end{cases}
$$
Taking the degree 0 part 
we get a projective complex of $A$ modules
with just the desired cohomology, establishing the result in this
case.

Next, suppose more generally that $M=S\otimes_A N$ is relatively projective.
As $A$ is concentrated in degree 0 the module $N$ splits
as a direct sum of modules of different degrees, so we may
as well assume that $N$ is concentrated in a single degree 
$d =\reg N = \reg \sF$.
 
We may identify the category of sheaves on $\Spec A$
with the category of graded $A$-modules concentrated in 
degree 0.
Since $\sF=\widetilde{M}$ we have 
$\sF=(\pi^* N(d))\otimes \sO_{\PP^n_A}(-d)$, whence
$$
R^i\pi_* \sF = N(d)\otimes R^i\pi_*\sO_{\PP^n_A}(-d).
$$
As before,
$R^i\pi_*\sF$ is nonzero for at most one $i$ (which is either
0 or $n$) so $R\pi_*\sF$ is represented by
this $A$-module,
as a complex concentrated in cohomological degree 0 or $n$.

On the other hand
$\FF$ is a complex of projective $A$-modules so, 
by Lemma \ref{hom-tensor},
$$
\Hom_S(\FF, S\otimes N)=\Hom_S(\FF, S(-d))\otimes N(d).
$$
By our hypotheses on $I$ the module $S/I$ is $A$-projective, and
 $\Hom_S(\FF, S(-d))$ is a resolution of $\Hom_A(A/I, A)(n+1-d)$.
This is split as a sequence of $A$-modules, so tensoring with $N(d)$
commutes with taking homology. It follows that the tensor product
has the desired cohomology.

Finally, we allow $M$ to be an arbitrary finitely generated graded 
$S$-module. By Hochschild's Theorem \ref{relative syzygy} $M$
has a finite $A$-split relatively projective resolution. We do induction
on the ``relative projective dimension'' of $M$,
that is, the minimal length $m$ of such a resolution. The case $m=0$ is
the case where $M$ is relatively projective. By Proposition \ref{elementary}
$M$ then has the form $S\otimes_AN$, and this is the case we
have already treated. 

Otherwise, let
$$
0\to M'\to S\otimes N\to M\to 0
$$
be the beginning of an $A$-split resolution of $M$ by relatively projective modules,
so that the relative projective dimension of $M'$ is $m-1$ and
the short exact sequence is $A$-split. Because each module
in $\FF$ is 
relatively projective, and because the previous sequence is
$A$-split, it follows that
$$
0\to \Hom(\FF, M')\to \Hom(\FF, (S\otimes N))\to \Hom(\FF, M)\to 0
$$
is an exact sequence of complexes, and thus forms a distinguished
triangle in the derived category. By induction, the 
two complexes on the left are the pushforward complexes of the 
appropriate sheaves, and it follows that the one on the right is too.
\qed

\section{Regularity and BGG}

In this section we make generalize
 Castelnuovo-Mumford regularity to the relative case.
Propositions \ref{base-change} and  \ref{flatness}
have no analogues in the absolute case, and seem to
be new properties of flat sheaves.
 
\medbreak
\noindent{\bf Relative Castelnuovo-Mumord Regularity}

 If $N$ is a graded
$S$-module such that $N_d=0$ for $d\gg 0$, then we define
the \emph{regularity} of $N$ by
$$
\reg(N) = \max\{d\mid N_d\neq 0\}.
$$
For a finitely generated graded $S$-module, on the other hand,
we define
$$
\reg(M) = \max_i ( \reg \Tor^S_i(A,M) - i).
$$
In the absolute case, when $A$ is a field,
this is the usual definition (see for example
Eisenbud [Eis04]). When $M$ is a finitely generated module
concentrated in finitely many degrees, so that both definitions
above are applicable, they coincide. This
can be proved by re-interpreting
regularity in terms of local
cohomology, as follows:

\begin{proposition}\label{regularity defs} 
$$
\reg(M) = \max_i ( \reg \Tor^S_i(A,M) - i) = \max_j ( \reg \HH^j_{(x)}M + j)
$$
and for each $i$
$$
\reg \Tor^S_i(A,M) - i \geq \reg \HH^{n+1-i}_{(x)}M +n+1-i
$$
where $(x)$ denotes the ideal $(x_0,\dots,x_n)\subset S$ and 
$\HH^j_{(x)}$ is the local cohomology functor. \hfill \qed
\end{proposition}

Since the proof is similar to that in
the classical case (see
for example Eisenbud [2004, Corollary 4.5]), we omit it.

\begin{example} {\bf Regularity is not semicontinuous.}\\
The following example is based on the fact that an ideal can
have higher regularity than a complete intersection with generators
of the same degrees. A convenient case is given by Caviglia [Cav04]:
If we set $i=(x_0^3,x_1^3,x_0x^2+x_1x_3^2)\subset S_0=K[x_0,\dots,x_3]$
the $S_0/i$ has regularity 7, whereas the ideal generated by a regular sequence
of 3 cubics has regularity 6. Consider the family of ideals, with
parameter $t$, given by
$$
I=(x_0^3,x_1^3,x_0x^2+x_1x_3^2+tx_2^3)\subset S=K[t][ x_0,\dots,x_3],
$$
The generators of $I$ form a regular sequence, and it follows that
the regularity of $S/I$  in our sense is 6, whereas the regularity of the
special fiber $S/(I,t)$ is 7. 
\end{example}

In general the regularity can go either up or down on specialization, though the
following result shows that in the flat case it can only go down.

\begin{proposition}\label{base-change}
In the situation above, suppose that $B$ is an $A$-algebra.
If either $M$ or $B$ is $A$-flat,
 then the regularity of $M\otimes_A B$ as a module
over $S\otimes_A B$ is at most the regularity of $M$ as an $S$-module.
If $A$ is local, $B$ is the residue class field of $A$, and $M$ is $A$-flat,
then $\reg_{S\otimes_AB} M\otimes_AB = \reg_S M$.
\end{proposition}

\begin{proof}
Let 
$$
\FF:\quad \cdots\to F_i \to\cdots \to F_0
$$ 
be a graded free resolution of $M$ as an $S$-module. 
By definition $M$ has 
regularity $r$
if and only if, for every $d> r$ and every $i\geq 0$,
the sequence 
$$
(F_{i+1}\otimes_SA)_{d+i} \to\cdots\to (F_{0}\otimes_SA)_{d+i} \to 0,
$$  
  obtained by taking the degree
$d+i$ component of the first $i+1$ steps of $\FF\otimes_SA$,
is exact up to and including the term $(F_{i}\otimes_SA)_{d+i} $.
The modules in this sequence are all $A$-free, so if the
sequence is exact, then it is split exact as a sequence of $A$-modules. 
This condition is
stable under base change to $B$: the sequence
$$
(F_{i+1}\otimes_SB)_{d+i} \to\cdots\to (F_{0}\otimes_SB)_{d+i} \to 0,\leqno{(*)}
$$  
is split exact as a sequence of $B$-modules.

Let $T=S\otimes_AB$.
If $M$ or $B$ is $A$-flat, then the sequence $\FF\otimes_ST\cong \FF\otimes_AB$
is exact, and thus it is a $T$-free resolution of $M\otimes_ST\cong M\otimes_AB$.
The regularity of $M\otimes_AB$ as a $T$-module is
thus computed from the exactness of the sequences
$(*)$, and the desired inequality follows.

Under the hypothesis of the last statement of the Proposition the other
inequality follows from Nakayama's Lemma.
\end{proof}

\goodbreak
\bigbreak
\noindent{\bf Bernstein-Gel'fand-Gel'fand Correspondence}

Using this notion of regularity, we easily generalize 
the BGG correspondence.
We need it in the form given for the case where $A$ is
a field in
which we studied in
[EFS, 2003] following [BGG, 1978]. The proofs
are essentially the same as in the case treated there.
For the reader's convenience we
formulate the necessary statements,
which involve a pair of adjoint
functors
$$
cplx(S) \stackrel{\rTo^R}{\lTo_L} cplx(E).
$$
 
To the graded $S$-module $M$ we associate 
$$
\RF(M):\quad \cdots \to \Hom_A(E,M_d) 
\to \Hom_A(E,M_{d+1}) \to \cdots,
$$
which is (in all interesting cases) an infinite complex of $E$-modules
with the differential  given by the action of $t=\sum_{j=0}^n
x_j\tensor e_j$. Since 
$$
\Hom_A(E, -)\cong E\otimes \wedge^{n+1} W\otimes -
$$
canonically, this complex is isomorphic up to shift in grading 
to the complex 
$$
\cdots \to E\otimes M_d 
\to E \otimes M_{d+1} \to \cdots
$$
of the  introduction. 
The definition extends  to the case where $M$ is
a bounded complex of $S$-modules by 
taking $\RF(M)$ to be the total
complex of the induced double complex. 

Similarly $\LF$ is defined by associating to a graded $E$-module 
$P=\oplus_j P_j$ the complex
$$
\LF(P):\quad \cdots \to S\tensor_A P_j \to 
S \tensor_A P_{j-1} \to \cdots.
$$ 
with $P_j$ concentrated in degree $j$, and $S\tensor_AP_j$
the term of homological degree $j$. For example the complex 
$\LF(E)$ has the form
$$
0 \to S \to S^{n+1}(1) \to \bigwedge^2(S^{n+1}(1)) 
\to\cdots \to \bigwedge^{n+1}(S^{n+1}(1))\cong S(n+1)\to 0,
$$ 
which is the Koszul complex.
\medskip

\begin{theorem}\label{regularity}  
Let $M$ be a finitely generated $S$-module.
\begin{itemize}
\item[a)] The truncated complex
$\RF(M_{\ge s})$ is acyclic for $s \geq \reg M$. 
\item[b)] Suppose $s\geq \reg M$, and set
$$
P= \ker(\Hom_A(E,M_{s+1}) \to
 \Hom_A(E, M_{s+2})).
$$ 
The complex 
$$\LF(P): \ldots \to S\tensor_A P_{s+2} \to S\otimes_A P_{s+1} \to 0 $$
is acyclic, and is a resolution of $M_{>s}$ via the map
$$
S\tensor_A P_{s+1} 
= S\tensor\Hom_A(\Lambda^{0} V,M_{s+1}) \to M_{> s}\hbox{ where } 
f\tensor_A \varphi \mapsto f\varphi(1).
$$
\end{itemize}\qed
\end{theorem}

We use Proposition \ref{base-change} to deduce a property
of flat sheaves and modules:

\begin{proposition}\label{flatness} If the sheafication 
$\sF=\widetilde M$ is flat over $\Spec A$, 
then for $s> \reg(M)$ the module
 $M_s$ is projective as an $A$-module. Moreover, the complex
$$
    E\otimes M_{s-1} \to E\otimes M_s \to E\otimes M_{s+1}
$$
is split exact at $E\otimes M_s$ as a complex of $A$-modules,
and thus 
$$
P^{s-1}:=\ker\bigl( E\otimes M_s \to E\otimes M_{s+1}\bigr)
$$ 
is projective as an $A$-module.
\end{proposition}

We will use the following well-known result, whose proof we provide
for the reader's convenience:

\begin{lemma}\label{splitness criterion}
Suppose that $A$ is a commutative Noetherian ring, and 
$$
F:   0\to M' \to M\to M'' \to 0
$$
is an exact sequence of finitely generated
$A$-modules, with $M$ projective. If
$F\otimes_A A/P$ is exact for all maximal ideals $P$ of $A$,
then $F$ is split exact.
\end{lemma}

\noindent It is easy to give examples where the conclusion fails if
we do not require $M$ to be projective.

\medbreak
 \noindent{\it Proof of Lemma \ref{splitness criterion}.\/} By the Noetherian and finite generation hypotheses, the 
formation of $\Ext^1(M'', M')$ commutes with localization.
Thus it is enough to prove the Lemma after localizing $A$ at $P$. We may factor
out the largest summand of $M$ contained in $M'$, and thus assume that $F$
is the beginning of a minimal free resolution of $M''$. 
In particular, we may assume that $M'$ is contained in $PM$. But now tensoring
with $A/P$, we see that $M' \otimes A/P = 0$, so $M' = 0$ by Nakayama's Lemma. 
\qed

\medbreak
 \noindent{\it Proof of Proposition \ref{flatness}.\/}From the exactness of the sequence
$$
0\to \HH^0_{(x)} M \to M \to \bigoplus_{d\in \ZZ} \pi_* \sF(d) 
\to \HH^1_{(x)} M \to 0
$$
we see that
$$
M_s\to \pi_* \sF(s)
$$ 
is an isomorphism  for $s> \reg(M)$. 

For the first statement of the Proposition it thus 
suffices to show that $\pi_*\sF(d)=M_d$
is flat, given that $\sG:=\sF(d)$ is flat. This holds very generally:
If $I\subset A$ is any ideal, then by the flatness of $\sG$ the map
the map 
$\pi^{-1}I\otimes_{\pi^{-1}A} \sG 
\to \pi^{-1}A\otimes_{\pi^{-1}A} \sG=\sG$
is 
a monomorphism, where we have written $A$ instead of $\sO_{\Spec A}$. 
Writing $\sO$ for $\sO_{\PP^n_A}$ we see that
the composite map
$$
\pi^* I \tensor \sG =
\pi^{-1} I\tensor_{\pi^{-1}A}\sO\otimes_\sO\sG=
\pi^{-1} I\tensor_{\pi^{-1}A} \sG 
\to 
\pi^{-1} A \tensor_{\pi^{-1} A} \sG = 
\sG
$$
is a monomorphism. Since the map $I\otimes_A \pi_*\sG \to \pi_*\sG$
is the direct image of $\pi^*I \otimes \sG \to \sG$,
it too is a monomorphism, so $\pi_*\sG$ is flat.

We next must show the sequence
$$
E\otimes M_{s-1} \to E\otimes M_s \to E\otimes M_{s+1}
$$
is split exact in the middle (that is, the kernel of the right hand map
is equal to the image of the left hand map, and is a direct summand
of $E\otimes M_s$.

Exactness is simply the statement of Theorem \ref{regularity} a). By Proposition \ref{base-change}, the regularity
of $M\otimes A/P$ is also $<s$, so the sequence
$$
E\otimes M_{s-1} \otimes A/P \to E\otimes M_s\otimes A/P \to E\otimes M_{s+1}\otimes A/P
$$
is also exact. By the first part of the Proposition,
 $M_s$ (and thus also $E\otimes M_s$) is $A$-projective.
 Lemma \ref{splitness criterion} now completes the proof.
\qed

\section{The Tate Resolution}

We next define and construct
Tate resolutions of a sheaf $\sF$ on $\PP^n_A$. 
When $A$ is a field, the Tate resolution is
defined in [EFS 2003] by joining the complex $\RF(M_{>s})$, 
for $s>\reg M$, to a free
resolution of $\ker(\Hom_A(E,M_{s+1})\to \Hom_A(E,M_{s+2}))$
to make a doubly infinite acyclic complex of free modules,
but in the relative case, when $\sF$ is not flat,
 the modules in $\RF(M_{>s})$ are
not free, and the situation is somewhat more subtle. 
In the general case, the Tate resolution
 will be 
a doubly infinite acyclic complex of graded 
relatively projective $E$-modules.  

If $T=\oplus_jE\otimes_AN_j$ is a graded, relatively projective
$E$-module, with $N_j$ concentrated in degree $j$, then
we write 
$$
T_{\preceq s} := T/(\bigoplus_{j>s}E\otimes_AN_j)
\cong  \bigoplus_{j \leq s} E\otimes_A N_j $$ 
with $N_j$ concentrated in degree $j$,
and call it the $s$-th \emph{generator truncated} quotient of $T$.
This notion extends to complexes: if 
$$
\TF:\quad \cdots \to T^{i-1} \to T^{i} \to\cdots
$$
is a complex of graded $E$-modules of the form $E\otimes_AN$ as above, we denote by
$\TF_{\preceq s}$ the complex formed from the $T^i_{\preceq s}$.
It is a quotient complex of $\TF$.

A visual idea of $\TF_{\preceq s}$ may help the reader through
what follows. We exhibit a sort of \emph{Betti diagram} for $\TF$, in the sense
of Bayer and Stillman (see for example [MACAULAY2]), by putting all the $N^i_j$ in a grid with
each column corresponding to a particular $T^i$:
\begin{diagram}[small]
&T^{s-1}&T^s&T^{s+1}\cr
&\vdots&\vdots&\vdots\cr
\cdots&{\bf N^{s-1}_{s-2}}&{\bf N^s_{s-1}}&{\bf N^{s+1}_s}&\cdots\cr
\cdots&{\bf N^{s-1}_{s-1}}&{\bf N^s_s}&N^{s+1}_{s+1}&\cdots\cr
\cdots&{\bf N^{s-1}_{s}}&N^s_{s+1}&N^{s+1}_{s+2}&\cdots\cr
&\vdots&\vdots&\vdots\cr
\end{diagram} 
Here the maps go from left to right. Those represented by matrices of 
scalars (degree 0 in $E$) go diagonally up and to the right, 
$N^i_j\to N^{i+1}_j$; the maps represented by matrices of linear
forms in the variables of $E$ (degree -1)
are horizontal, $N^i_j\to V\otimes N^{i+1}_{j+1}\subset E\otimes N^{i+1}_{j+1}$, 
etc. We have
indicated the part of the complex belonging to $\TF_{\preceq s}$ 
by putting it in boldface.

As usual, we work with a finitely generated graded $S$-module $M$ and
the coherent sheaf $\sF=\widetilde M$ on $\PP^n_A$ associated to it.
We set 
$$
P^s=\ker (E\otimes M_{s+1}\to E\otimes M_{s+2}).
$$ 
Note that $P^s$ differs from the module 
$\ker(\Hom_A(E,M_{s+1})\to \Hom_A(E,M_{s+2}))$ used
in the last section only by a shift of degrees, since
$\Hom_A(E,M_d)\cong \Hom_A(E,A)\otimes M_d\cong
E\otimes \wedge^{n+1}W\otimes M_d$, 
and $\wedge^{n+1}W\cong A(n+1)$.

We define a \emph{Tate resolution} of $\sF$ to be an
acyclic, doubly infinite
complex $\TF=\TF(\sF)$ of graded relatively projective
$E$-modules $T^i=E\otimes_A N^i$ such that
$\TF_{\preceq s}$ is a resolution of $P^s$  when 
$s\gg 0$. A \emph{projective Tate resolution} is a Tate
resolution in which all the $T^i$ are projective as $E$-modules or, 
equivalently, all the $N^i$ are projective as $A$-modules.
Here  $P^s$ is concentrated in 
degree $s+1$, and in cohomological degree $s$; thus
we require $\TF_{\preceq s}$ to be concentrated in 
homological degrees $\leq s$. In our construction, 
we will produce a Tate resolution such that
$\TF_{\preceq s}$ is a resolution of $P^s$ for all 
$s\geq \reg M$ and all the $N^i$ are finitely generated.
We will also produce a projective Tate resolutions with
the first property; but
it may happen that there are no projective Tate resolutions
where all $N^i$ are  finitely generated. 

\medbreak
\noindent{\bf Construction of Tate Resolutions}

Let $M$ be any finitely generated graded $S$-module
such that $\widetilde M = \sF$. Choose an integer $s\geq \reg M$,
and let $P^s=\ker (E\otimes M_s\to E\otimes M_{s+1})$ as above.
Choose a relatively projective resolution 
$$
\cdots \to T^{s-1}\to T^s\to P^s\to 0.
$$
 Let $\TF$ be the complex
$$
\TF:\qquad \cdots \to T^{s-1}\to T^s\to E\otimes_AM_s\to E\otimes_AM_{s+1}\to \cdots
$$
\begin{corollary}\label{construction in flat case}
The complex $\TF$ above is a Tate resolution for $\sF$. If $\sF$ is
flat over $A$ and the $T^i$ are chosen to be projective over
$S$, then the complex $\TF$ is a projective Tate resolution that
is $A$-split.
\end{corollary}

\begin{proof}
Since $s\geq \reg M$,
Theorem \ref{regularity} shows that the right hand
part of this complex is acyclic, and the left hand part is acyclic
by construction. For $u\geq s$
the truncations $\TF_{\preceq u}$ are
resolutions of the desired form. If $\sF$ is flat then the right-hand 
part 
$$
\RR(M_{\geq s}): E\otimes_A M_s\to E\otimes_AM_{s+1}\to\cdots
$$
 is an $A$-split exact projective complex by Proposition \ref{flatness},
In particular $P^s$ is $A$-projective, so the left-hand part
$$
\cdots\to T^{s-1}\to T^s
$$
is also $A$-split.
\end{proof}

When $\sF$ is not flat over $A$ we will construct a projective
Tate resolution by choosing a projective resolution $\TF$ of $P^s$, for some
$s\geq \reg M$, and then extending it inductively to
a Tate resolution $\TF(\sF)$ by adding generators in degrees
higher than $s$  so that $\TF = \TF(\sF)_{\preceq s}$.

\begin{proposition}\label{extending T} Suppose that
$s\geq \reg M$ and let
 $\TF$ be a projective resolution of $P^s$.
Let $F: \ldots \to F^{s-1}\to F^{s} \to F^{s+1} \to M_{s+1}(s+1) \to 0$
be a projective resolution as of $M_{s+1}(s+1)$ as an $A$-module.
The mapping cone $\TF'$ of
$$ 
\TF \to E \tensor_A F(-s-1)
$$
is a projective resolution of $P^{s+1}$, and 
$\TF=\TF'_{\preceq s}$. For $j\leq s$ the inclusion $\TF\subset \TF'$
induces an equality
$(A\otimes_E\TF')_j=(A\otimes_E\TF)_j$.
                      
If $A$ is local and both $ \TF$ and $F$
are  minimal free resolutions, then so is $\TF'$.
\end{proposition}

\begin{proof} 
By Theorem \ref{regularity} the sequence
$$
0 \to P^s \to E\tensor M_{s+1} \to P^{s+1} \to 0
$$
is exact.
The complex 
$E\tensor F(-s-1)$
is a resolution of the $E$-module $\to E\tensor M_{s+1}$. Thus the given mapping
cone is  a
projective resolution of $P^{s+1}$ and differs frm $\TF'$  by the
addition of free modules generated in degree exactly $s+1$.

In the local case, if $\TF'$
and $F$  are  minimal, then the free modules in
$\TF'$ are all generated in degrees $\leq s$,
the generator degree of $P^s$, since the elements of $E$ 
have degree $\leq 0$. Thus the comparison maps from 
$\TF(\sF)_{\preceq s}$ to $E\tensor F$ cannot contain
degree 0 components, and the mapping cone is again minimal.
\end{proof}

We can now complete the construction of projective Tate resolutions in
general:
Let $M$ be a finitely generated graded $S$-module, and let $s>r:=\reg M$
be an integer.
Let $\sF=\widetilde M$ be the corresponding coherent sheaf on  $\PP^n_A$.
Set $P^s=\ker E\otimes M_{s+1}\to E\otimes M_{s+2}$ as above, and 
suppose that for $s\geq r$ the $A$-module $M_s$ has projective dimension
at most $d$ (we allow the case $d=\infty$).
Let $\TF^r$ be a projective resolution of $P^r$. Supposing that
$\TF^s$ has been constructed for $s\geq r$, let $\TF^{s+1}$ be the
result of extending $\TF^s$, as in Proposition \ref{extending T}, 
using a projective resolution of length at most $d$.
The union $\TF=\TF(\sF)$ of the $\TF^s$ is a Tate resolution of $\sF$.

If $d<\infty$ then each $T^u$ will be finitely generated, but 
otherwise this may not be the case.

\medbreak
\noindent{\bf Small Tate resolutions}

It is interesting, especially from a computational point of view,
to construct  Tate resolutions that are as small as possible.
For this we need a basic fact about the existence of
resolutions of the form
$$
 (*)\qquad \cdots \to \bigoplus_j E\otimes_AN^{-1}_j
\to \bigoplus_j E\otimes_AN^0_j
\to P \to 0.
$$

\begin{proposition}\label{koszul regularity}
Let $P$ be a graded $E$-module. If $s\leq t$
are integers such that
$P_j=0$ unless $s\leq j\leq t$, then: 
\begin{itemize}
\item[(a)] The module $P$ has a (not necessarily projective)
resolution of the form $(*)$ such that
$$
N^{i}_j=0\hbox{ unless }i\leq 0 \hbox{ and }s+i\leq j\leq t+i.
$$
\item[(b)]
If $P$, as an $A$-module, has a projective resolution of length $d$
then $P$, as an $E$-module, has a projective resolution of the form
$(*)$ such that
$$
N^{i}_j=0\hbox{ unless }i\leq 0 \hbox{ and }s+i\leq j\leq \min(t,t+d+i).
$$
\end{itemize}
\end{proposition}

It is easy to give examples where these inequalities are sharp.

\begin{proof} 

First consider the case $s=t$. The Cartan-Eilenberg
resolution of $A$ over $E$,
$$
{\bf CE}: \quad\cdots \to E\otimes (\Sym_2 W)^*\to E\otimes W^*\to E \to A\to 0
$$
has $(\Sym_m W)^*$ is concentrated in degree $-m$
for each $m$. The modules in the complex ${\bf CE}\otimes_AP$ 
have the form $E\otimes N^i_{t+i}$, for $i\leq 0$, where
$N^i_{t+i}$ is a direct sum of finitely many copies of $P$.
We have 
$$
H^i({\bf CE}\otimes_AP) =\Tor^E_{-i}(A,P)=0\hbox{ for }i<0
$$
so ${\bf CE}\otimes_AP$ is a resolution satisfying part $(a)$.

If 
$$
{\bf F}:\cdots \to F^{-1}\to F^0\to P\to 0
$$
 is a projective resolution of $P=P_t$, so that 
each $F^i$ is concentrated in degree $t$, then the 
homology of the
projective complex
$H^i({\bf CE}\otimes_A {\bf F})$ is also equal to $\Tor_{-i}^E(A,P)=0$,
so ${\bf CE}\otimes_A {\bf F}$ is a projective resolution of
$A\otimes P$ as an $E$-module. If $F$ has length
$d$, the $i$-th term is
$$
T^{i}=\bigoplus_{\substack{u,v\leq 0\\u+v=i\\  -d\leq v}}
{CE^u\otimes_AF^v} = 
\bigoplus_{\substack{u,v\leq 0\\u+v=i\\  -d\leq v}}
E\otimes_A (\Sym_{-u} W)^*\otimes_AF^{v}.
$$
It has generating projective $A$-module 
$$
N^{i}_{t+j}=(\Sym_{-j} W)^*\otimes_AF^{i-j},
$$
which is nonzero only when $-j\geq 0$, $-d\leq i-j\leq 0$---that is,
in the range $j=t+i,\dots,\min(t, t+d+i)$, as required
for part $(b)$.

For the general case, we do induction on the difference
$t-s$, so we may assume, for either statement $(a)$ or $(b)$,
that we have complexes 
$T'$ and $T''$ of the desired form for
the modules $P'=P_s$ and $P''=P/P_s$. We can construct
from these a resolution $T$ of $P$ whose $i$-th term
is $(T')^i\oplus (T'')^i$, and this complex will
have the desired form.

To construct $T$ in  case $(a)$, where the terms of 
$T''$ are not assumed projective,
but simply of the form $E\otimes_A N$,
we must use the fact that, as an $A$-module, $P_s$ is a summand
of $P$ so that
the exact sequence 
$$
0\to P'\to P\to P''\to 0
$$
is $A$-split.

We now construct
the resolution of $P$ inductively: given 
\begin{diagram}[small]
&&E\otimes_AN'&& &&E\otimes_AN''\cr
&&\dTo&& &&\dTo\cr
0&\rTo &Q'&\rTo &Q&\rTo &Q''&\rTo& 0
\end{diagram}
with the bottom sequence $A$-split.
Because $E\otimes N$ is relatively
projective we may lift the map $E\otimes N\to Q''$ to a map $E\otimes N\to Q$,
and use it and the given map on the left to construct a diagram
\begin{diagram}[small]
0&\rTo &K'&\rTo &K&\rTo &K''&\rTo& 0\cr
&&\dTo&&\dTo &&\dTo\cr
0&\rTo&E\otimes_AN'&\rTo&E\otimes_A(N'\oplus N'') &\rTo&E\otimes_AN''&\rTo&0\cr
&&\dTo&&\dTo &&\dTo\cr
0&\rTo &Q'&\rTo &Q&\rTo &Q''&\rTo& 0
\end{diagram}
where the modules $K',K,K''$ are by definition the kernels of the 
vertical maps below them. The middle row and the bottom row
are compatibly $A$-split, so the exact sequence 
$0\to K'\to K\to K''\to 0$ is also $A$-split, and we 
may continue the construction.

The construction of $T$ in the case where
the complexes $T'$ and $T''$ are projective
works the same way, and is of course standard.
 \end{proof}

\begin{corollary}\label{range} 
Let $\sF$ be a coherent sheaf on $\PP^n_A$ represented by 
a finitely generated 
graded $S$-module $M$ of regularity $r$. The sheaf
$\sF$ has a Tate 
resolution $\TF$ with terms $T^i=\bigoplus_jE\otimes_AN^i_j$
such that 
$$
N^i_j \hbox{ is nonzero only in the range }i-n \le j \le i, 
$$
while if $i> r$ then 
$$
N^i_j \hbox{ is nonzero only if }j=i.
$$

If  every
$A$-module has projective dimension $\leq d\leq \infty$,
then there is a Tate resolution of $\sF$
such that each $N^i_j$ is projective, and 
$$
N^i_j \hbox{ is nonzero only in the range }i-n \le j \le \min(r, i+d), 
$$
while if $i> r$ then 
$$
N^i_j \hbox{ is nonzero only in the range }i \le j \le i+d.
$$
\end{corollary}

In the general case, where
$A$ does not have finite global dimension
and $\sF$ is not flat, the modules
$\TF(\sF)^i$ in the projective Tate resolution
we have constructed will not be finitely generated.
However,
the module $(\TF(\sF)_{\preceq s})^k$ is
finitely generated. To compute
$(A\otimes_E \TF(\sF))_u$, for any $u$, it suffices to
compute finitely many terms of
 the complex $\TF^s=\TF(\sF)_{\preceq s}$, for example
as a resolution of $P^s$, for any $s\geq\max(\reg M, u)$.

\begin{proof} 
The limits  follow immediately from
our construction, using Proposition \ref{koszul regularity} and
Proposition \ref{extending T}.
\end{proof}

We can express the restrictions of
Corollary \ref{range} pictorially in terms of a 
Betti diagram: 
$$
\begin{tabular}{rrr rrl l lll} 
$\cdots \quad \cdots $&$ \cdots $&$ \cdots$&$ 
\cdots $&$  N_0^n\quad\ldots $&
$ N^{r-1}_{r-1-n} $&
$$&$$&$$\cr\vspace{4pt} 
$\cdots\quad\cdots$&$  \cdots   $&$ \cdots $&$  
\antiddot$&$\cdots \quad \cdots $&
$ \quad\vdots$&
&&$$\cr\vspace{4pt}
$\cdots\quad\cdots $&$\cdots $&$ N^{1}_{0}$&$ \cdots  $&$\cdots\quad \cdots $&
$N^{r-1}_{r-2}$
&$ $&$ $&$  $ \cr\vspace{4pt}
$\cdots\quad\cdots $&$N^{0}_{0} $&$ \cdots$&$ \cdots $&$\cdots\quad \cdots $&
$N^{r-1}_{r-1} \;\, N^{r}_{r}$&
$ \cdots$&$ \cdots $&$\hspace{-3pt} N^{s}_{s}$ \cr\vspace{4pt}
$\cdots\;N^{-1}_{0}$&$\cdots $&$ \cdots$&$ \cdots $&$\cdots \quad \cdots$&$ 
N^{r-1}_{r} \;\cdots $
&$\cdots$&$\hspace{-3pt} N^{s-1}_{s} $&$ $\cr\vspace{4pt} 
$\antiddot\quad\cdots$&$ \cdots $&$\cdots $&$ \cdots $&$\cdots\; N^{r-2}_{r}
$&$\;\cdots \quad \cdots $
&$\antiddot $&&$ $\cr \vspace{4pt}  
$\cdots \quad\cdots $&$ \cdots $&$ \cdots $&$\cdots $&$\antiddot \quad\cdots  $&
$\;\cdots\quad\antiddot$
&$$&$$&$$\cr
\end{tabular} 
.$$
In case $A$ is local, the $N^i_j$ are free, and we sometimes put just the
numbers $\rank N^i_j$ into the table. 

The $i$-th column in the table corresponds to the term $\TF(\sF)^i$,
Maps pointing directly to the right
are linear in the variables $e_i$, while maps of higher degree in
the $e_i$ point to the right and downwards.
The degree 0 maps in $\TF$ go
up and to the right, along the $45^0$ diagonals; thus the complex
$(A\otimes_E \TF_{\preceq s})_d$ is
$$
\begin{matrix} 
\hspace{3.3cm}\nearrow \hspace{0cm}\cr
\hspace{2.75cm}N^{1}_d \hspace{0.55cm}  \cr
\hspace{2.2cm}\nearrow\hspace{1.1cm} \cr
\hspace{1.65cm}N^0_d\hspace{1.65cm}\cr 
\hspace{1.1cm}\nearrow\hspace{2.2cm} \cr
\hspace{0.55cm}
N^{-1}_d \hspace{2.75cm} \cr
\hspace{0cm}\nearrow \hspace{3.3cm}\cr
\end{matrix}
$$
The empty spaces in the table, and in particular
the absence of dots above the top row shown, 
indicate that the terms that might have occured there
are 0; these are the restrictions of Corollary \ref{range} (in case
$d=\infty$).

\section{The Beilinson Monad}

In this section we use a relatively projective Tate resolution $\TF$ for
$\sF$ to construct a relative Beilinson monad for
$\sF$: it is a complex $\UF(\sF)$ 
of $\pi_*$-acyclic coherent sheaves on $\PP^n_A$,
 whose only homology is $\sF$.
The construction is parallel to that given in 
[ES 2003], but involves new subtleties, especially
when $\sF$ is not flat over $A$.

Let $U$ denote the universal subbundle  on $\PP^n_A$.
$$
0 \to U \to W\otimes \sO \to \sO(1) \to 0,
$$
so that $U = \Omega^1_{\PP^n_A/\Spec A}(1)$.
We will define an additive functor
$$\UF: \{\hbox{relatively projective } E\hbox{-modules}\} \to coh(\PP^n_A)$$
from  the category of relatively projective $E$-modules
to the category of sheaves on $\PP^n_A$ sending
$E\otimes N_j$ to $\wedge^jU\otimes \pi^*N_j(j)$. 
Here we have written $\pi^*N_j(j)$ instead of $\pi^*N_j$
because we identify the category of 
quasicoherent sheaves on $\Spec A$ with the category of
$A$-modules concentrated in degree 0.

To define the action of $\UF$ on morphisms, note that
$$
\Hom_E(E\tensor N_j, E\tensor N_\ell)= 
\Lambda^{\ell-j}V\tensor_A \Hom_A(N_j(j),N_\ell(\ell)), $$
where as always it is understood that $N_j$ sits in degree $j$ and $N_\ell$ in
degree $\ell$.
A morphism $\phi\in \Lambda^{\ell-j} V\tensor \Hom_A(N_j(j),N_\ell(\ell))$
defines, for each $p$ a morphism 
$\wedge^{p}W\otimes \pi^*N_j(j)\to \wedge^{p+j-\ell}W\otimes \pi^*N_\ell(\ell)$,
and since these morphisms commute with the differentials of the 
Koszul complex, they induce a map $\UF(\phi)$ that makes the
following diagram commute (see [ES03] for more
details):
$$
\xymatrix{ E\tensor_A N_j \ar[d]_\phi &  & \Lambda^{-j} U \tensor \pi^*N_j(j)\ar[d]_{\UF(\phi)} 
\ar[r]& \Lambda^{-j} W
  \tensor \pi^*N_j(j) \ar[d]_{\phi} \cr
E\tensor_A N_\ell  &  & \Lambda^{-\ell} U \tensor \pi^*N_\ell(\ell)\ar[r] & \Lambda^{-\ell} W
  \tensor \pi^*N_{\ell}(\ell) \cr} 
$$

Note that
$\Lambda^{-j} U = 0$ unless $0 \ge j \ge -n$.
Because of this the functor $\UF$ will annihilate
many terms of a Tate resolution
$\TF(\sF)$
so
$$
\UF(\sF):= \UF(\TF(\sF))
$$
is a bounded-above complex
that  depends only on $n+1$ diagonals of $\TF(\sF)$.
In particular $\UF(\sF)=\UF(\TF(\sF)_{\preceq s})$ for any 
$s\ge \max(0,\reg(M))$.

Recall that a \emph{monad} for $\sF$ is a complex $U$
of sheaves
whose homology is 0 except for $H^0(U)=\sF$.

\begin{theorem} 
Let $M= \sum_d M_d$ be a finitely generated graded
$S$-module and let $\sF$ be the
  associated coherent sheaf on $\PP^n_A$. If
$\TF(\sF)$ is a Tate resolution for $\sF$, then
the complex
$\UF=\UF(\TF(\sF))$
is a bounded-above  monad for $\sF$ whose
terms are $\pi_*$-acyclic
coherent sheaves. There are relatively projective
Tate resolutions $\TF$ for $\sF$ such that 
$\UF$ is a finite complex. On the other hand,
if $A$ is local
 then $\UF$ may be determined by $\sF$, uniquely up to isomorphism,
by the additional requirement that 
that $\TF(\sF)_{\preceq s}$ is a
minimal free resolution for every $s>\reg M$.
Furthermore, if 
\begin{enumerate} 
\item $\sF$ is $A$-flat, or
\item $\pd A < \infty$
\end{enumerate}
holds, then we can choose a projective
Tate resolution such that
the complex $\UF(\sF)$ is finite.
\end{theorem}

\begin{proof} The sheaf $\wedge^jU$ is resolved by
a truncated Koszul complex,
$$
0\to  \Lambda^{n+1} W \tensor \sO(i-n-1) \to 
\cdots \to \Lambda^{i+1} W \tensor \sO(-1) \to
\Lambda^i U \to 0.
$$
We deduce, just as in the case where $A$ is a field,
that $R^j\pi_* \Lambda^i U =0$ for $i>0$ and all $j$
and that 
$R\pi_* \Lambda^0 U =
R^0\pi_* \sO_{\PP^n_A}=A$.
Thus all terms in the complex $\UF(\sF)_{\preceq s}$ are
 $\pi_*$-acyclic.   

To prove that $\HH^*(\UF(\sF))=\HH^0(\UF(\sF)) \cong \sF$
we use a double complex argument. 
If $G$ is a graded $S$-module or a complex of such, we let
$\widetilde{G}$ stand for the sheafication of a (complex) of 
graded $S$-modules $G$.

Theorem \ref{regularity} says that
$\widetilde \LF(P^s\tensor\Lambda^{n+1}W)$ 
is a complex with $\sF$ as its only homology. 
Consider the subcomplex $(\TF(\sF)_{\preceq s})_{\ge -n}$ 
formed by replacing each 
$$
(T^i)_{\preceq s}=\bigoplus_{j\leq s}E\otimes_A N^i_j
$$
by the submodule 
$$
((T^i)_{\preceq s})_{<-n}
=\bigoplus_{j\leq s, \ell\geq n+j}\wedge^\ell V\otimes_A N^i_j
$$
of 
elements in  $\TF(\sF)_{\preceq s}$ of
degree $\ge -n$. We apply the functor $\widetilde \LF$
to this subcomplex.
The resulting double complex has as vertical complexes
$\LF(((T^i)_{\preceq s})_{<-n})$
that are sums of sheafifications  of  possibly truncated Koszul complexes
$\widetilde  L((E\tensor N^i_j)_{\ge -n})$
as in the following diagram:
$$ 
\begin{tabular}{ccc}
$0 $&& 0 \cr\vspace{2pt}
$\uparrow$&&$\uparrow$\cr\vspace{2pt}
\small{$ \sum\limits_{\ell-j=n}  \pi^*(\Lambda^{\ell} V \tensor
N^i_j(j))\tensor \sO(n)$}&$\rightarrow $&\small{$
\sum\limits_{\ell-j=n}  \pi^*(\Lambda^{\ell} V \tensor
N^{i+1}_j(j))\tensor \sO(n)$}\cr\vspace{2pt}
$\uparrow$&&$\uparrow$\cr\vspace{2pt}
\small{$\sum\limits_{\ell-j=n-1}  \pi^*(\Lambda^{\ell} V \tensor
N^i_j(j))\tensor \sO(n-1)$}&$\rightarrow $&\small{$
\sum\limits_{\ell-j=n-1}  \pi^*(\Lambda^{\ell} V \tensor
N^{i+1}_j(j))\tensor \sO(n-1)$}\cr \vspace{2pt}
$\uparrow$&&$\uparrow$\cr\vspace{2pt} 
\end{tabular}
$$
The  complex $\widetilde  L((E\tensor N^i_j)_{\ge -n})$
a resolution of
$\Lambda^{-j} U \tensor \pi^*(N^i_j(j)) \tensor \sO(n+1)$
in case $0 \ge j \ge -n$, and otherwise it is exact.
Thus the vertical homology of the double complex is the complex 
$\UF(\sF)\tensor \sO(n+1)$.

On the other hand, since
$\TF(\sF)_{\preceq s+1}$
is a resolution,
the horizontal homology 
of
$\LF(((T^i)_{\preceq s})_{<-n})$
 is the complex $\widetilde L(P^s)$.
In the notation of Theorem \ref{regularity}, 
we have $P^s=P(-n-1)$, so by that result
the homology of $\widetilde L(P^s)$ 
 is
$\sF(n+1)$.
 
A diagram chase in the double complex proves
$\HH^*(\UF(\sF)\tensor \sO(n+1)) = \HH^0(\UF(\sF)\tensor \sO(n+1)) \cong
\sF(n+1)$ as desired.

The last statement follows from
 Corollary \ref{range} and Proposition \ref{flatness}.
\end{proof}

\begin{corollary} 
If $\TF(\sF)$ is any Tate resolution for $\sF$
then the complex
$\pi_*\UF(\TF(\sF))$ represents $R\pi_*\sF$ in the derived category
of bounded-above, finitely generated complexes of $A$-modules.
In particular  $\pi_*\UF(\sF)$ has no homology in negative degrees
and $R^i\pi_*\sF =\HH^i \pi_* \UF(\sF)$ for $i\ge 0.$
In case $A$ is local and we choose $\TF(\sF)_{\preceq s}$
to be a minimal free resolution, then 
$\pi_*\UF(\sF)$
is the unique minimal free representative.
\end{corollary}

\begin{proof} We have exact complexes
$$
\ldots \to \UF(\sF)^{-2} \to  \UF(\sF)^{-1} \to \sB \to 0,
$$
$$ 
0 \to \sK \to \UF(\sF)^{0} \to \UF(\sF)^{1} \to \ldots
$$
and
$$
0 \to \sB \to \sK \to \sF \to 0.
$$
Since $\UF(\sF)$ is $\pi_*$-acyclic and $R^{n+1}\pi_* \sG=0$ for any coherent
sheaf we get an exact complex
$\ldots \to \pi_*\UF(\sF)^{-2} \to \pi_*\UF(\sF)^{-1} \to \pi_* \sB \to 0$
and $R^i\pi_* \sB = 0$ for $i>0$.
Thus the result follows from 
$$0 \to \pi_*\sB \to R\pi_* \sK \to R\pi_*\sF \to 0$$
and $R\pi_* \sK = \pi_* \UF(\sF)^{\ge 0}$. 
\end{proof}

{\it Completion of the proof of Theorem \ref{main}}. 
Since
$$\pi_* \Lambda^\ell U =\begin{cases} 0, &\hbox{ if }\ell\not=0\cr
A, &\hbox{ if } \ell=0
\end{cases}$$
we obtain
$$R\pi_*\sF=\pi_*\UF(\sF)=((\TF(\sF)_{\preceq s} \tensor_E A)_0,$$ 
as desired. 
\hfill \qed

\begin{remark} Let $\sB^0=\coker(\UF(\sF)^{-1}\to \UF(\sF)^{0})$. The proof 
shows that
$$0\to \pi_* \sB^0 \to \pi_* \UF(\sF)^1 \to \ldots \to
  \pi_*\UF(\sF)^n \to 0$$
is a bounded complex representing $R\pi_* \sF$.
\end{remark}

Grothendieck's motivation for introducing $R\pi_* \sF$ as a complex was to
make a base change property true. It is amusing to note that the property
follows directly from our construction.

\begin{corollary}[Base change]\label{base change} Suppose $\sF$ is flat over $A$. Then
  $R\pi_*\sF$
commutes with base change in the sense that $\varphi^*R\pi_*(\sF)$ represents 
$R\pi_*(\widetilde \varphi^* \sF)$
for any ring homomorphism $\varphi \colon A \to A'$ and the induced diagram
$$
\xymatrix{
\Spec A' \times \PP^n \ar[r]^{\widetilde\varphi} \ar[d]
& \Spec A \times \PP^n \ar[d] \cr
\Spec A' \ar[r]^{\varphi} & \Spec A
}
$$
\end{corollary}

\begin{proof} In the flat case Corollary \ref{construction in flat case} 
produces a Tate resolution that commutes with arbitrary base change,
since it is split exact as a complex of $A$-modules.
\end{proof}

\begin{corollary}\label{shape}
Suppose $A$ is local with maximal ideal $\gm$ and $\sF$ flat over $A$. Then
the $k$-th summands in the minimal Tate resolution and
Beilinson monad  are
$$
(\TF(\sF))^k= \sum_{i=0}^n E\tensor A(-k+i))^{h^i(k-i)},
$$ 
and
$$
(\UF(\sF))^k= \sum_{i=0}^n (\wedge^{i-k}U)^{h^i(k-i)},
$$ 
where $h^i(k-i)= \dim_{A/\gm} \HH^i(\sF(k-i)\tensor_A A/\gm)$,
as in the case where $A$ is a field.
\end{corollary}

\section{Examples}

\begin{example}
\label{example1}
To exhibit our technique in a simple and natural case,
 we take the versal deformation $\sF$ of the bundle
$\sE=\sO\oplus \sO(-2)$ where $\sO=\sO_{\PP^1}$ is the structure sheaf
of the projective line over a field $K$, and compute the
complex $R\pi_*\sF$. 

The base space of this deformation
has as tangent space
$$
\Ext^1(\sE, \sE)\cong \HH^1(\sHom(\sE,\sE)) =\HH^1(\sHom(\sO,\sO(-2)) =K,
$$
and since the deformations are unobstructed the base space of the
versal deformation is the germ of $\AAA^1=\Spec K[a]$. We thus work over 
$A=K[[a]]$.

By Corollary \ref{shape} the Betti diagram of the Tate resolution 
$\TF(\sF)$ is 
$$
\begin{tabular}{r|cccccc}
$i-j\setminus i$&-2&-1&\,0&\,1&\,2&\,3\cr\hline
1&6&4&2&1 \cr
0& & &1&2&4&6
\end{tabular}
$$
Hence the minimal representative of $R\pi_*\sF$ is the complex
$$
\xymatrix{0\ar[r]& A^1\ar[r]^\alpha&  A^1\ar[r]&  0\cr
}
$$
for some map $\alpha$ (which is not hard to guess, but which 
we will compute to illustrate our method in this easy example.)

From the Betti diagram we see directly that the regularity of $\sF$ 
is 2, and we can compute $R\pi_*\sF$ starting from the map
of  modules 
$$
\phi_2: E\otimes (\HH^0\sF(2))\to E\otimes (\HH^0\sF(3))
$$
over the exterior algebra $E$.

We write $x,y\in W=\pi_*\sO(1)$ for fiber coordinates on $\PP_A^1$,
where now $\sO$ denotes the structure sheaf of $\PP^1_A$,
and $e,f$ for their dual coordinates in 
$E$.
The sheaf $\sF(2)$ is an extension
$$
0\to \sO \to \sF(2)\to \sO(2)\to 0,
$$
Lifting a basis for $\HH^0\sO(2)$,
we may choose a basis of the free $A$-module $\HH^0\sF(2)$ 
denoted by $1, x^2, xy, y^2$. In terms of this basis, a presentation
matrix may be written
$$\bordermatrix{
&\cr
1&-ax&0\cr
x^2&y&0\cr
xy&-x&y\cr
y^2&0&-x\cr
}
$$
We choose as basis of $\HH^0\sF(3)$ the elements 
$$
x\cdot 1,\ y\cdot 1,\
x\cdot x^2,\ x\cdot xy,\ x\cdot y^2,\ y\cdot y^2.
$$
From the given relations we see that $y\cdot xy=x\cdot y^2.$
However, $y\cdot x^2 = x\cdot xy+a(x\cdot 1).$ Thus, in
terms of these bases, the 
map $\phi_2$ has matrix
$$
\bordermatrix{
&1&x^2&xy&y^2\cr
x\cdot 1&e&af&0&0\cr
y\cdot 1&f&0&0&0\cr
x\cdot x^2&0&e&0&0\cr
x\cdot xy&0&f&e&0\cr
x\cdot y^2&0&0&f&e\cr
y\cdot y^2&0&0&0&f
}
$$
The further syzygy matrices are
$$
\phi_1=\begin{pmatrix}
ef&af&0\cr
0 &e&0\cr
0 &f&e\cr 
0 &0&f\cr  
\end{pmatrix},
\phi_0=\begin{pmatrix}
a&e&f\cr
e&0&0\cr
f&0&0
\end{pmatrix},
\phi_{-1}=\begin{pmatrix}
fe&0&0&0\cr
af&e&f&0\cr
0 &0&e&f\cr 
\end{pmatrix},\quad
$$
Hence $\phi_{-1}={}^t\phi_1$ (transposed in the sense
appropriate to the exterior algebra), and
$\phi_{-2}={}^t\phi_2$ as well.

Finally, by Theorem \ref{main}, $R\pi_*\sF=(A\otimes_E\TF(\sF))_0$ is the complex
$$
\xymatrix{0\ar[r]& A^1\ar[r]^{(a)}&  A^1\ar[r]&  0.\cr
}
$$
Similar computations can be done very quickly for much larger examples
by \Mac . 
\end{example}

\begin{example}[\bf Vector Bundles on $\PP^1$]

\label{example1a}

More generally, consider the family of globally generated
vector bundles rank $r$ and degree $d$ on $\PP^1$.
The most special bundle in this family is
$\sF_0= \sO^{r-1}  \oplus \sO(d)$. Every other bundle in this family arises
as an extension
$$0 \to \sO^{r-1} \to \sF_a \to \sO(d) \to 0,$$
with $a \in \Ext^1(\sO(d),\sO^r)\cong \HH^1(\sO^{r-1}(-d))\cong
\HH^0(\sO(d-2)^{r-1})^*$.
Thus as a base space of this family we may choose $\Spec A$ with
$$A=K[a_i^s, \,0 \le i \le d-2, 1 \le s \le r-1],$$
and then take $\sF$ to be the universal extension on $\PP^1\times A$.

We will specify $\sF$ explicitly via its Beilinson monad. 
The Tate resolution has Betti diagram 
$$
\begin{tabular}{r|ccc c cccc}
$i-j\setminus i$&-d-1     &-d     &-d+1       & \ldots &-2      &-1&0&1\cr\hline
1             &d(r-1)+r &d(r-1) &(d-1)(r-1) & \ldots & 2(r-1) & r-1 &     &      \cr
0             &         &1      &2          & \ldots &  d-1   & d   & d+r & d+2r \cr
\end{tabular}$$
at the special point $0 \in \Ext^1(\sO(d),\sO^{r-1})$. A few examples computed
with \Mac \ make it possible to guess the pattern of the differentials, which
we now going to verify.
The differentials of 
$$\xymatrix{\TF(\sO):& \ldots \ar[r]^{{}^tC^3}&  E^2(3) \ar[r]^{{}^tC^2}& E(2) \ar[r]^{(ef)} & E
\ar[r]^{C^1}& E^2(-1) \ar[r]^{C^2} &E^3(-2) \ar[r]^{C^3} & \ldots} $$ 
are given by special $(\ell+1) \times \ell$ \emph{Toeplitz}
(or \emph{Hankel}) matrices
$$C^\ell =\begin{pmatrix}
    e & 0 &        & \ldots & 0 \cr
    f & e & \ddots &        &\vdots \cr
    0 & f & \ddots  & \ddots &  \cr
\vdots & \ddots & \ddots & e & 0 \cr
 & &\ddots & f & e \cr
0 & \ldots& & 0 & f\cr  
\end{pmatrix}
$$ 
and their transposed. Hence $\TF(\sF)$ is a deformation of the complex $\TF(\sF_0)=\oplus_1^{r-1} \TF(\sO) \oplus
\TF(\sO)[d](d)$  build from
Toeplitz matrices. 
To describe the most relevant piece we consider
for pairs $(k,\ell)$ with $k+\ell= d$  the $k\times \ell$ Hankel matrices
$$B_{k\ell}^s :=\begin{pmatrix} 
a^s_0 & a^s_1 & a^s_2 & \ldots & a^s_{\ell-1} \cr
a^s_1 & a^s_2 & a^s_3 & \ldots & a^s_{\ell} \cr
\vdots & \vdots & \vdots & \ddots & \vdots \cr
a^s_{k-1}&a^s_k & a^s_{k+1}& \ldots& a^s_{k+\ell-2} \cr 
\end{pmatrix}.
$$

\begin{proposition}\label{T on P1} The $-(k+1)$ th differential
$$\xymatrix{ E(k+2)^{(k+1)(r-1)}\oplus E(k+1)^{\ell} \ar[r]^{\quad D} &
  E(k+1)^{k(r-1)}\oplus E(k)^{\ell+1}},
$$ 
in the Tate resolution of $\sF$ for $1 \le k \le d-1$ and $\ell=d-k$ is given by the block matrix 
$$D=D^{-k-1}=
\left( \begin{tabular}{cccc|c}
${}^t C^{k} $&$ 0 $&$ \ldots $&$ 0 $&$ B_{k\ell}^1 $\cr\vspace{1pt}
$0 $&$ {}^t C^{k} $&$ \ddots $&$ \vdots $&$B_{k\ell}^2 $\cr\vspace{1pt}
$\vdots $&$ \ddots$&$\ddots $&$ 0 $&$ \vdots $\cr\vspace{1pt}
$0 $&$ \ldots $&$ 0 $&$ {}^t C^{k} $&$ B_{k\ell}^{r-1} $\cr\hline
$0$&$ \ldots $&$ 0 $&$ 0 $&$ -C^\ell $\cr
\end{tabular}\right)
$$
(in suitable coordinates). 
\end{proposition}

\begin{proof} We first prove that the matrices
  $D^{-d},D^{-d+1},\ldots,D^{-2}$ define a complex $\TF$. Indeed
$D^{-k} \cdot D^{-k-1}=0$ holds because
$$ {}^tC^{k-1}\cdot B_{k\ell}^s - B_{k-1,\ell+1}^s\cdot C^\ell =0.$$
With the relative Beilinson monads we can recover the corresponding 
coherent sheaf from any two consecutive matrices: 
The monad $\UF(\TF(-k)[-k])$ is the total complex of a double complex
$$\xymatrix{
\sO(-1)^{k(r-1)} \ar[r] & \sO^{(k-1)(r-1)} \cr
\sO(-1)^\ell \ar[u]^{B_\ell}\ar[r] & \sO^{\ell+1} \ar[u]_{B_{\ell+1}} \cr}$$
whose rows do not depend on the parameters. Here
$$B_\ell = \begin{pmatrix}
B^1_{k\ell}\cr \vdots \cr B^{r-1}_{k\ell}
 \end{pmatrix}.$$  
The homology of the top row, which is a subcomplex, is
$$\sO(-k)^{r-1}=\ker(\sO(-1)^{k(r-1)} \to \sO^{(k-1)(r-1)}),$$ while the
coorresponding quotient complex, the bottom row, has homology
$$\coker(\sO(-1)^\ell \to \sO^{\ell+1}) = \sO(\ell)=\sO(d-k).$$
 Thus $\HH^* \UF(\TF(-k)[-k])=\HH^{0}\UF(\TF(-k)[-k])=: \sF_k$ and the
 homology of
$\UF(\TF(-k)[-k])$ fits into a short exact sequence

\noindent
$$
(*_k)\hspace{2cm}0 \to \sO(-k)^{r-1}\longrightarrow \sF_k \longrightarrow \sO(d-k)
\to 0 \,.\hspace{2cm}
$$
\noindent
Since the complex defined by the $D^{-d},\ldots,D^{-2}$ has the right Betti
number
and the matrices have linearly independent rows, we conclude that they form part of
the Tate resolution of a single sheaf $\sF$, and, that $\sF_k=\sF(-k)$ with the extensions
$(*_{k+1})=(*_k)\tensor \sO(-1)$. To prove that $\sF$ is the universal
extension, it suffices to prove this for anyone of the sheaves
$\sF_k=\sF(-k)$. We choose $\sF_d=\ker(\sO(-1)^{d(r-1)}\oplus \sO \buildrel
D^{-d} \over \longrightarrow 
   \sO^{(d-1)(r-1)})$. 
The boundary map in
$$ \to \HH^0(\PP^1_A, \sF_d) \to \HH^0(\PP^1_A,\sO) \to
\HH^1(\PP^1_A,\sO(-d)^{r-1}) \to $$   
is the composition  
$$\HH^0(\PP^1_A,\sO){\buildrel B_1\over\longrightarrow}
\HH^0(\PP^1_A,\sO^{(d-1)(r-1)})\cong \HH^1(\PP^1_A,\sO(-d)).$$
So the boundary map vanishes at a point $a$ iff all coordinates $a^s_i$
vanish at $a$. We conclude, that  the $a^s_i$ represent linearly independent extension
classes, and, since we have the right number $(d-1)(r-1)$ of parameters, that
$(*_d)$ 
is the universal extension.  
\end{proof}  

\begin{corollary} Let $\sF$ on $\PP^1_A$ be the universal extension 
$$
0\to \sO^{r-1} \to \sF \to \sO(d) \to 0
$$
Then for each $k$ in the range  $1\le k \le d-1$ the direct image complex of
$\sF(-k-1)$ is
$$
R\pi_* \sF(-k-1):  0 \to A^{d-k}\; {\buildrel B_{d-k}\over \longrightarrow} \;
  A^{k(r-1)} \to 0 
\hbox{ with } B_{d-k}=\begin{pmatrix} B^1_{k,d-k} \cr
  \vdots\cr B^{r-1}_{k,d-k} \cr \end{pmatrix}$$
Outside this range the direct image complexes are concentrated in one degree.
\end{corollary}

\end{example}

\begin{example}[\bf Strata in the case $(d,r)=(6,3)$]
The corollary allows to describe the loci of extension classes of a given
splitting type in $\Ext^1(\sO(d),\sO^{r-1})$ by rank conditions on the matrices
$B_{k}$ in the various direct image complexes.

We treat the example $(d,r)=(6,3)$. The possible splitting type correspond to
partition of $d$ into at most $r$ parts. In our special case this are the
following strata with an arrow $p \to q$ indicating that the strata $p$ lies
in the closure of the strata $q$:

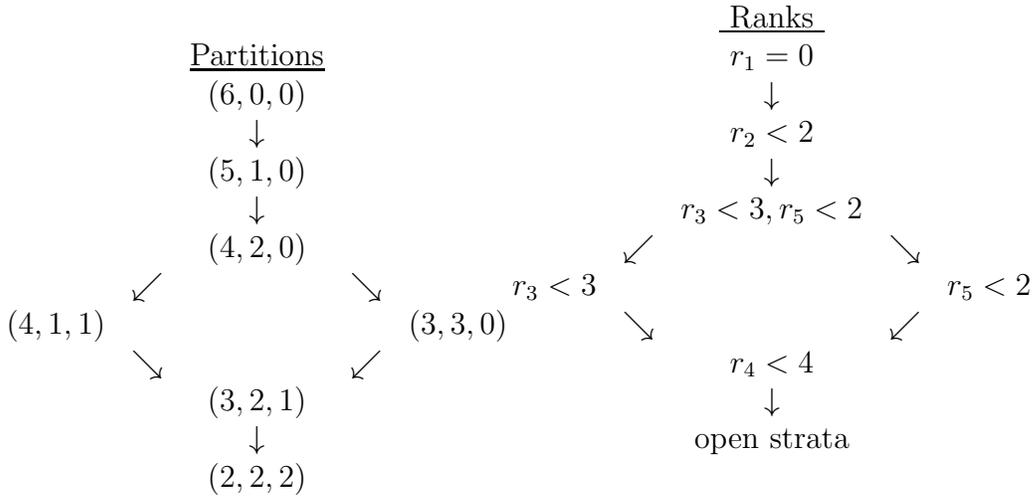
\noindent
\vspace{-0.5cm}
$$\begin{matrix}
&&&&\cr&&&&\cr
&&\hbox{\underline{Partitions}}&&\cr 
&& (6,0,0) && \cr
&&  \downarrow && \cr   
&& (5,1,0) && \cr
&&  \downarrow && \cr 
&& (4,2,0) && \cr
&\swarrow&& \searrow & \cr 
(4,1,1)&&&& (3,3,0)\cr
&\searrow &&\swarrow&  \cr
&& (3,2,1) && \cr
&&  \downarrow && \cr
&& (2,2,2) && \cr
\end{matrix}
\begin{matrix}
&&\hbox{ \underline{ Ranks } }&& \cr
&& r_1=0 && \cr
&&  \downarrow && \cr
&& r_2<2&& \cr
&&  \downarrow && \cr
&& r_3<3 , r_5< 2 && \cr
&\swarrow&& \searrow & \cr
\,r_3<3&&&& r_5<2\cr
&\searrow &&\swarrow&  \cr
&& r_4<4 && \cr
&&  \downarrow && \cr
&& \hbox{open strata}&& \cr
\end{matrix}
$$

\noindent
In which strata an extension $a$ lies is determined by the ranks $r_i = \rank
B_{i}(a)$ of the $ (d-i)(r-1)\times i$ matrices $B_{i}$ evaluated at $a$.
$$
\begin{tabular}{r|cccccccc}
$i-j\setminus i$&-7&-6&-5&-4&-3&-2&-1&0\cr\hline
1&15&12&10&8&6&4&2& \cr
0&  &1&2&3&4&5&6&9\cr
\end{tabular}
$$
Indeed by the Base Change Theorem \ref{base change} and its Corollary
\ref{shape} the $r_i$ 
determine the dimensions $h^0(\PP^1,\sF_a(-i-1))$, which in turn determine the
splitting type according to the following elementary Lemma.

\begin{lemma} Let $\sE$ be a vector bundle on $\PP^1$, and let
$$h=h_\sE: \ZZ \to \ZZ,\, n\mapsto h^0\sE(n)$$
be its Hilbert function. Then
$$\sE \cong \oplus_{j \in \ZZ} \sO(-j)^{h''(j)}\,,$$
where $h''(j)=h(j)-2h(j-1)+h(j-2)$ denotes the second difference function.
\hfill \qed\end{lemma}

The claim on the strata indicated in the table above follows.
Note that there is no single matrix, whose rank determine all splitting types,
and for one strata there is not a single matrix on which a rank condition
gives the defining equations.

To exhibit the beautiful pattern in this family of matrices more visibly, we 
drop the upper index notation $a^s_i$  and use coordinates
$a_0,\ldots,a_4,b_0,\ldots,b_4$ instead. With this notation we have
$$
B_5=\begin{pmatrix} 
a_0 & a_1 & a_2 & a_3 & a_4 \cr\hline
b_0 & b_1 & b_2 & b_3 & b_4 \cr
\end{pmatrix}
,
B_4=\begin{pmatrix}
a_0 & a_1 & a_2 & a_3  \cr
a_1 & a_2 & a_3 & a_4 \cr\hline
b_0 & b_1 & b_2 & b_3  \cr
b_1 & b_2 & b_3 & b_4 \cr
\end{pmatrix}
, B_3=\begin{pmatrix}
a_0 & a_1 & a_2  \cr
a_1 & a_2 & a_3  \cr
a_2 & a_3 & a_4  \cr\hline
b_0 & b_1 & b_2  \cr
b_1 & b_2 & b_3 \cr
b_2 & b_3 & b_4 \cr
\end{pmatrix}.
$$
There are various relations between the ideals of minors of these matrices. The
most interesting one is the primary decomposition of ideal $3\times 3$ minors
of the square matrix $B_{4}$:
$$minors(3,B_{4})=minors(3,B_{3})\cap minors(2,B_{5}).$$
We discovered this relation by computation using Macaulay2, which provides
a proof in a few positive characteristics; the relation was recently
proven non-computationally by Moty Katzman [Ka05, Section 3].

In terms of projective geometry the closed strata have the following descriptions as
cones over projective varieties:
$$\begin{matrix}&&\hbox{\underline{Geometry:}}&&&\hbox{codim}&\hbox{degree}\cr
&& \emptyset &&&10 & 1 \cr
&&  \downarrow &&&& \cr
&& S(4,4) &&& 7 & 8\cr
&&  \downarrow && \cr
&& Sec(\PP^1 \hookrightarrow \PP^4)\times \PP^1 &&&5 & 12 \cr
&\swarrow&& \searrow &&& \cr
Sec(S(4,4))&&&& \PP^4\times \PP^1&4& 15, 5\cr
&\searrow &&\swarrow&  \cr
&& Sec^3(S(4,4)) && &1&4\cr
&&  \downarrow && \cr
&& \PP^9&& &0&1\cr
\end{matrix}
$$
Here $S(4,4)\subset \PP^4\times \PP^1 \subset \PP^9$ denotes the 2-dimensional rational normal scroll
defined by the $2\times 2$ minors of the matrix
$${}^tB_{2}= 
\left(\begin{tabular}{cccc|cccc}
$a_0 $&$ a_1 $&$ a_2 $&$ a_3 $&$ b_0 $&$ b_1 $&$ b_2 $&$ b_3 $\cr  
$a_1 $&$ a_2 $&$ a_3 $&$ a_4 $&$ b_1 $&$ b_2 $&$ b_3 $&$ b_4 $\cr  
\end{tabular}\right),
$$
and $Sec(X)$ respectively $Sec^3(X)$ refers to the secant respectively
3-secant variety of $X$. Note that the fibers of $S(4,4) \to \PP^1$ are
rational normal curves of degree $4$, and that 
$$S(4,4)\cong \PP^1\times \PP^1
{\buildrel |(4,1)|\over\hookrightarrow} \PP^9.$$

\end{example}

We return to the general case of the versal deformation of 
$\sO^{r-1}\oplus\sO(d)$.
The analysis above shows in general that, in the deformations of 
$O(d)\oplus \sO^{r-1}$, the stratum of bundles isomorphic to
$\sO(d-1)\oplus \sO(1)\oplus \sO^{r-2}$ is isomorphic to the 
rational normal scroll 
$$
S(d-2,\dots, d-2)\cong (\PP^1 \hookrightarrow \PP^{d-2})\times \PP^{r-2}
\subset \PP^{dr-d-r}.
$$
We would like to have a geometric description of the strata in general, for example in
terms of secant constructions such as 
$$
Sec^b\bigl(Sec^a (\PP^1 \hookrightarrow \PP^{d-2})\times \PP^{r-2}\bigr).
$$

Though we don't have such a geometric description, we can at least give
an algebraic one. 
If a closed point $p\in \Spec A$ corresponds to a bundle with a certain splitting type,
then from the cohomology of the bundle we see that the matrices in the 
pushforward of the $R\pi_*\sF(\ell)$, for various $\ell$, satisfy some rank 
conditions at $p$. We conjecture that the determinantal equations derived
from these rank conditions actually generate the radical ideal of closure
of the stratum of vector bundles that split in this way. The following stronger statement
includes the case of certain unions of strata:

\begin{conjecture} Let $B_1,\ldots,B_{d-1}$ be the non-trivial matrices of the direct
image complexes of the versal deformation of $\sO^{r-1}\oplus\sO(d)$.
For any collection of positive integers $r_{k_1},\ldots,r_{k_s}$, the ideal
$$
\sum_{t=1}^{s} minors(r_{k_t},B_{k_t})
$$
is radical.
\end{conjecture}
The case $r=2$ was asserted by Room [1938], and proven
in characteristic 0 by Peskine and Szpiro. See 
Conca [1998] for a general proof and other references. 

The minimal primes of the ideal
$$
\sum_{k=1}^{d-1} minors(r_k,B_k)
$$
are easy to describe, and are (radicals of) ideals of the same form. 
First of all, the locus
of extensions
$$
0\to \sO^{r-1}\to \sE\to \sO(d) \to 0
$$
on $\PP^1_K$ such that $\sE$ has a given splitting type,
$\sE\cong \oplus \sO(a_i)$,
is always irreducible. Its closure is locus
of extensions such that each twist of $\sE$ 
has at least as many global sections as the corresponding
twist of $\oplus \sO(a_i)$. The corresponding prime
ideal is thus the radical of
the corresponding sum of ideals
$minors(r_k,B_k)$. Conversely, any sum of
the ideals $minors(r_k,B_k)$ defines the locus of extensions such that 
 various twists of $\sE$ have at least a certain
number of independent global sections. With some care
one can give the irredundant decompositions in terms of splitting types.

In particular, the prime ideals
of the form
$$
\rad\bigl(\sum_{k=1}^{d-1} minors(r_k,B_k)\bigr)
$$
are precisely the primes that define the closures
of the strata of points $p\in \Spec A$ where 
$\sF_p$ has a given 
splitting type. If the conjecture is true, of course, these
sums of determinantal ideals are already radical.

One can see from this analysis that an ideal $minors(r_{k_t},B_{k_t})$
can have components of different
dimensions; in particular, it need not be Cohen-Macaulay.

It is worth remarking that one can also treat this family of
examples without exterior methods. For example, from the exact sequence
$0\to \sO^{r-1}(-k)\to \sF(-k) \to \sO(d-k)\to 0$ we get a triangle
of direct images that expresses 
$R\pi_*\sF(-k)$ as the mapping cone of a certain map
$$
R\pi_*\sO(d-k) \to R\pi_*\sO^{r-1}(-k)[1].
$$ At most one
of the modules $R^i\pi_*\sO(d-k)=\HH^i(\sO(d-k))$ is nonzero, and 
similarly for $\sO^{r-1}(-k)$, so each of
$
R\pi_*\sO(d-k)$
and 
$ R\pi_*\sO^{r-1}(-k)[1]
$
 reduces to a single free module. 
In the
``interesting'' range $-d\leq k\leq -2$ where both these
modules are nonzero,
the map between the modules is the connecting homomorphism
$\HH^0(\sO(d-k))\to  \HH^1(\sO^{r-1}(-k))$
that we have called $B_{d-k+1}$. This connecting homomorphism
is easy to compute concretely, especially since the computation
reduces to the case
$r=1$.

\begin{example}[\bf Blow-up of an elliptic singularity]
For an example that seems much harder to treat by
simple methods, consider the singularity defined by
$$ B=\{abc+a^4+b^4+c^4=0\} \subset \AAA^3.$$
The singularities of $B$ are resolved by blowing up the origin 
$\sigma:  \widetilde \AAA^3 \subset \PP^2\times \AAA^3\to \AAA^3$ once. 
We consider the strict $\overline B=\overline{\sigma^{-1}(B\setminus\{0\})} \subset \PP^2\times B$ and the total transform 
$B' =\sigma^{-1}(B)$ of $B$. The Tate resolution of $\sO_{\overline B}$
has Betti diagram
$$
\begin{tabular}{r|cccc c cccc}
$i-j\setminus i$
  &-4&-3&-2&-1&\,0&\,1&\,2&\,3&\,4\cr\hline
1 &15&12& 9& 6& 3 & 1 &   &   &   \cr
0 &27&21&15& 9& 4 & 3 & 6 & 9 & 12\cr
-1&21&15& 9&4 & 4 & 9 & 15& 21& * \cr
-2&15& 9& 4& 4& 9 & 15& 21& * & * \cr
\vdots&$\antiddot$&$\antiddot$&$\antiddot$&$\antiddot$&$\antiddot$&$\antiddot$&$\antiddot$&$\antiddot$&\cr
\end{tabular}
$$
with eventually periodic diagonals by [Eisenbud, 1980].

The Tate resolution of $\sO_{B'}$ looks quite different:
$$
\begin{tabular}{r|cccc c cccc}
$i-j\setminus i$
  &-4&-3&-2&-1&\,0&\,1&\,2&\,3&\,4\cr\hline
2 &10& 6&3&1& & & &  &  \cr
1 &15&8 &3& & && &  &  \cr
0 & 6&3 &1&1&1&3&6&10&15\cr
-1&  &  & & &3&8&15&24 &*\cr
-2&  &  & &1&3&6&10&* &*\cr
\end{tabular}
$$
with bounded $\TF^k$,
although $A=\sO_{B,0}$ has not finite projective dimension.
A closer inspection of the complexes gives
$$R\pi_* \sO_{B'}(k) = R\sigma_* \sO_{\widetilde \AAA^3}(k) \tensor_{\sO_{\AAA^3}} \sO_B,$$
a formula which holds, although the Base Change Theorem 
\ref{base change} does not apply.
So for $k\ge 0$ we have 
$$R^0\pi_* \sO_{B'}(k) = \gm_{\AAA^3,0}^k \tensor_{\sO_{\AAA^3,0}} \sO_B,$$
while
$$R^0\pi_* \sO_{\overline B}(k) =\gm_{B,0}^k.$$
\end{example}

\begin{example}[\bf Variety of Complexes]

One might hope that the direct image of a vector bundle
on $\PP^n_A$, say in the case of a local ring $A$, would have special properties compared
to an arbitrary complex of free $A$-modules. But it turns out that
such images are general; in fact one can get any complex as the push-forward
of quite a simple bundle:

\begin{theorem}\label{all occur} 
Every bounded minimal free complex
$$ 0 \to A^{\beta_0} \to  A^{\beta_1} \to \ldots \to A^{\beta_n} \to 0$$
over a local Noetherian ring arises as the direct image complex of a locally
free sheaf on $\PP^n_A$.
\end{theorem}

By flat base change, it suffices to prove the result for 
the \emph{generic complex}
$$
\FF: \qquad 0\to B^{\beta_0}\to \cdots\to B^{\beta_n}\to 0
$$
defined over the ring
$$B=\ZZ[a^p_{ij}; 1\le i \le \beta_{p+1}, 1\le j \le \beta_p,p=0,\ldots,n-1]
/(\, (a^{p+1}_{ij})(a^p_{jk})=0 \,), 
$$
with the map $B^{\beta_p}\to  B^{\beta_{p+1}}$ given by the map
with matrix $(a^p_{jk}).$
The ring $B$ is the affine coordinate ring of the
``variety of complexes''(see for example
[DS 1981]). The next Theorem is thus a strengthening of
Threorem \ref{all occur}.

\begin{theorem}\label{all occur 2}
The generic complex $\FF$ over the ring $B$ above is the
direct image of the versal deformation of the vector-bundle 
$$
\bigoplus_0^n (\wedge^p \Omega_{\PP^n_\ZZ/\ZZ})^{\beta_p}.
$$
\end{theorem}

\noindent Note that Example \ref{example1}
is
the special case of Theorem \ref{all occur 2} a complex of length 1 with free modules
of rank 1, 
$$
0\to B^1\to B^1\to 0.
$$

\begin{proof} To simplify the notation, we write
$\Omega^i$ for $\wedge^i\Omega_{\PP^n_\ZZ/\ZZ}$.
We first show that the ring $B$ (or more properly its completion
at the origin) is
 the base of the versal deformation $\sF$ of the vector bundle
$\sF_0=\bigoplus_{p=0}^n (\Omega^p)^{\beta_p}.$


\begin{proposition}
$$\Ext^1(\Omega^p,\Omega^q) = \begin{cases} \HH^1(\Omega^1)=\ZZ & \hbox{if
} 1\le
  q=p+1\le n \cr
0 & \qquad \hbox{ otherwise} \end{cases}$$
$$\Ext^2(\Omega^p,\Omega^q) = \begin{cases} \HH^2(\Omega^2)=\ZZ & \hbox{if
} 2\le
  q=p+2\le n \cr
0 & \qquad \hbox{ otherwise} \end{cases}$$
Thus the tangent space of the versal deformation of $\sF_0$
is
$$\Ext^1(\sF_0,\sF_0)=\bigoplus_{p=0}^{n-1}
\Hom(\ZZ^{\beta_p},\ZZ^{\beta_{p+1}}),$$
and the obstruction space is
$$\Ext^2(\sF_0,\sF_0)=\bigoplus_{p=0}^{n-2}
\Hom(\ZZ^{\beta_p},\ZZ^{\beta_{p+2}}).$$
The obstruction map is the given by composition
$$(\phi_0,\ldots,\phi_{n-1})\mapsto (\phi_1\cdot\phi_0,\ldots,
\phi_{n-1}\cdot\phi_{n-2}),$$
and the base space $B$ of the deformation is the coordinate
ring
of the variety of complexes of free modules of ranks $\beta_0,\ldots,\beta_n$.
\end{proposition}

\begin{proof} The generator of $\Ext^1(\Omega^{p-1},\Omega^{p})$
corresponds to the extension
$$
0 \to \Omega^{p}  \to \bigwedge^{p}(\sO^{n+1}(-1)) \to \Omega^{p-1} \to 0
$$
that appears in the Koszul complex, and the generator of 
$\Ext^2(\Omega^{p-1},\Omega^{p+1})$
can also be realized in the  Koszul complex as the extension
$$
0 \to \Omega^{p+1}  \to \bigwedge^{p+1}(\sO^{n+1}(-1)) 
\to \bigwedge^{p}(\sO^{n+1}(-1)) \to \Omega^{p-1} \to 0,
$$
which is thus the Yoneda product of the generators of 
$\Ext^1(\Omega^{p-1},\Omega^{p})$ and $\Ext^1(\Omega^{p},\Omega^{p+1})$.

The quadratic obstruction map
$$
\Ext^1(\sF_0,\sF_0)\to \Ext^2(\sF_0, \sF_0)
$$
is the map given by squaring, using the Yoneda product.
The only nonzero contributions come from the maps
$$\begin{matrix}
\Ext^1((\Omega^{p-1})^{\beta_{p-1}},(\Omega^{p})^{\beta_{p}})\times
\Ext^1((\Omega^{p})^{\beta_{p}},(\Omega^{p+1})^{\beta_{p+1}} )\cr
\downarrow \cr 
 \Ext^2((\Omega^{p-1})^{\beta_{p-1}},(\Omega^{p+1})^{\beta_{p+1}})
\end{matrix}
$$
which, with natural choice of bases, 
is matrix multiplication by the computation above, so we see that the obstructions
to second-order deformation are as claimed.

To see
that there are no higher terms in the equations of the base space 
of the versal deformation, we
observe
that $\sF_0$ can be
graded
by giving a summand $\Omega^p$ degree $p$. With this grading, each nonzero 
element of
$\Ext^1(\sF_0,\sF_0)$
is of degree 1 and the elements of $\Ext^2(\sF_0,\sF_0)$ are similarly of degree
2. Since the versal deformation must be homogeneous for this grading,
no higher-degree terms can occur in the equations.
\end{proof}

We also make use of the computation of the Tate resolution
associated to the sheaf $\Omega^p$: from [EFS 03], Proposition 5.5,
we know it has the form
$$\begin{tabular}{r|cccc c cccc}
$i-j\setminus i$
  &$\ldots$&$p-2$&$p-1$&$p$&$p+1$&$p+2$&\ldots&\cr
\hline
n &$\ldots$&$u_{p-2}$&$u_{p-1}$&$0$&$0$&$0$&$\ldots $&&            \cr
$\vdots$&       &      &        &          &&        &       &     & \cr
$p+2$ &           &        &         &$0$&           &        &      &      &  \cr
$p+1$ &           &        &         &$0$&  &        &      &      &  \cr
$p$ &           &        &         &$1$&           &        &      &      &  \cr
$p-1$ &           &        &       &$0$&           &        &      &      &  \cr
$p-2$ &           &        &         &$0$&           &        &      &      &  \cr
$\vdots$&       &      &        &          &&        &       &     & \cr
0 &$\ldots$&$0$&$0$&$0$&$v_{p+1}$&$v_{p+2}$&$\ldots$\cr
\end{tabular}
$$
with $u_{p-1}={n+1 \choose n+1-p}$ and $v_{p+1}={n+1 \choose p+1}$.

It remains to show that
the generic complex over $B$ is
$
R\pi_* \sF,$
where $\sF$ on $\PP^n_A$ is the versal deformation of 
$\sF_0=\oplus_{p=0}^n \oplus_{j=1}^{\beta_p} \Omega^p$ on $\PP^n_\ZZ$. 

From the Betti diagram above we see that the Tate resolution of $\sF_0$ has a Betti diagram of 
the following shape, where the entries not shown in the center of the table are all zero:
$$\begin{tabular}{r|cccc c cccc}
$i-j\setminus i$
  &-1&\,0&\,1&\,2&$\ldots$&n-1&\,n&n+1&$\ldots$\cr\hline
n &$\delta_{-1}$&$\delta_0$&$\delta_1$&$\delta_2$&$\ldots $&$\delta_{n-1}$&$ \beta_n$&           &\cr
n-1&       &      &        &          &          &$\beta_{n-1}$&        &            & \cr
$\vdots$&       &      &        &          &$\antiddot$&        &       &     & \cr
2 &           &        &         &$\beta_2 $&           &        &      &      &  \cr
1 &           &        &$\beta_1 $&        &           &        &       &     &    \cr
0 &           &$\beta_0 $&$\gamma_1$&$\gamma_2$&$\ldots$&$\gamma_{n-1} $&$\gamma_n$&$\gamma_{n+1}$&$\ldots$\cr
\end{tabular}
$$
Because the base ring $B$ is naturally graded,
we can use the result for local rings in 
Corollary \ref{shape}
to conclude that the 
Tate resolution of $\sF$  has the same entries.
We want to prove that the component
$$
E^{\beta_p} \to E^{\beta_{p+1}}
$$
of the differential $\TF(\sF)^p \to \TF(\sF)^{p+1}$,
which appears on the diagonal,
is given by the matrix $(a_{ij}^p)$.

Consider  the subspace 
corresponding to the factor ring $B_p= B/I_p$
with ideal $I_p$ generated by the linear forms
$(a_{ij}^q; q\not=p)$. The  restriction of the family $\sF$ to this subspace is
$$
\sF\tensor \sO_{B_p}\cong \bigl(\oplus_{q\not=p,p+1}
(\Omega^q)^{\beta_q}\bigr)
\oplus\sG
$$
where $\sG$ is the versal deformation of 
$$
\sG_0=(\Omega^p)^{\beta_p}\oplus (\Omega^{p+1})^{\beta_{p+1}}.
$$
The Tate resolution of $\sG$ is a deformation of the Tate resolution of 
$\sG_0$; the shape of its Tate resolution can be deduced from that
of $\Omega^p$ and that of $\Omega^{p+1}$ as above. 

We focus on the two differentials
$$\begin{matrix} E^{\beta_{p+1}(n-p){n+2 \choose n+1-p}}(n+1-p) \cr
\oplus \cr
E^{\beta_p{n+1 \choose n+1-p}}(n+1-p)
\end{matrix}
\buildrel   c \over \longrightarrow
\begin{matrix}  E^{\beta_{p+1}{n+1 \choose n-p}}(n-p)\cr \oplus \cr
 E^{\beta_p} \end{matrix} 
\buildrel d  \over \longrightarrow
\begin{matrix}
 E^{\beta_{p+1}}\cr \oplus \cr
E^{\beta_p{n+1 \choose p+1}}(-p-1)
\end{matrix},$$
with
$$c=\begin{pmatrix} c_1 & c_{12} \cr 0 & c_2\end{pmatrix} 
\hbox{ and }
d=\begin{pmatrix} d_1 & d_{12} \cr 0 & d_2\end{pmatrix}.$$
By [EFS, 2003] Proposition 5.5 we know that
$c_2$ is a direct sum of $\beta_{p}$ copies of the $1\times {n+1 \choose n+1-p}$
matrix consisting of all monomials of degree $n+1-p$ in the exterior 
variables $e_0, \ldots e_n$, and that $d_2$ consist of a direct sum $\beta_p$
copies of the ${n+1 \choose p+1}\times 1$ matrix consisting of all monomials
of degree $p+1$ in the exterior variables. (So $d_2\cdot c_2=0$ because
the composition has degree $n+2$.) Similarly $d_1$ is a direct sum of
$\beta_{p+1}$ copies of the 
$1\times {n+1 \choose n-p} $ matrix of all monomials of degree $n+1-(p+1)$
in $e_0,\ldots,e_n$ and $c_1$ consists of $\beta_{p+1}$ copies of the
linear syzygies matrix of these monomials.

The deformation sits in the components $c_{12}$ and $d_{12}$ of the matrices
corresponding to the extensions. We want to prove that 
$$d_{12}=(a_{ij}^p).$$  
If we take this choice for $d_{12}$ then we can build a complex by
 taking $c_{12}$ as a suitable matrix of bihomogeneous forms in the variables
$a_{ij}^p$ and $e_0,\ldots,e_n$, because 
$(e_0,\ldots,e_n)^{n+1-p} \subset (e_0,\ldots,e_n)^{n-p}$. The two differentials,
deformed in this way, extend to a map of (doubly infinite) resolutions,
and thus to a 
deformation of the whole Tate resolution. This defines a sheaf $\sG'$
on $\PP^n \times B_p$. 

We argue now directly that $\sG'$ over $B_p$ defined in this way is 
the versal deformation of $\sG_0$. Indeed for any other deformation 
$\sG''$ of $\sG_0$ over some base space $\Spec T$ 
the direct image complex of $\sG''$ induces a morphism 
$\varphi:\Spec T \to B_p$ 
by taking the 
substitution $\ZZ[a_{ij}^p] \to T$ obtained from the matrix $d_{12}''$
in the complex $R\pi_* \sG''$ on $\Spec T$, and
 $\sG'' \cong (\id_{\PP^n}\times \varphi)^* \sG$.
This proves the (semi) universal property of $\sG'$. It is universal
because the grading of $\sG'$ and its base ring $\ZZ[a_{ij}^p]$
given by degree in the $a_{ij}^p$
prevents there being any automorphisms except for conjugation of
the maps $(a_{ij}^p)$ by invertible matrices in the obvious way.
Thus $\sG'=\sG$.

The Base Change  Theorem  \ref{base change}
and the 
grading of $\sF$ and its base ring
$\ZZ[a_{ij}^p;p=0,\ldots,n-1]$ by degree in the $a_{ij}^p$ shows
that $(a_{ij}^p)$ occurs as a differential in
$R\pi_* \sF$, since
 $\sG$ is a summand of $\sF\tensor \sO_{B_p}$, as required.
\end{proof}
\end{example}

\begin{conjecture} Any bounded complex of finitely generated $A$-modules
$$
0\to B_0\to\cdots\to B_n\to 0
$$
arises as the direct image complex of a family of sheaves on $\PP^n_A$, which
can be taken  to be a deformation of the sheaf $\oplus_p (B_p\otimes\wedge^p \Omega_{\PP^n_A/A})$.
\end{conjecture}
The methods above suffice to prove this for a two-term complex.

\bigskip

\vbox{\noindent Author Addresses:\par
\smallskip
\noindent{David Eisenbud}\par
\noindent{Department of Mathematics, University of California, Berkeley,
Berkeley CA 94720}\par
\noindent{eisenbud@math.berkeley.edu}\par
\smallskip
\noindent{Frank-Olaf Schreyer}\par
\noindent{Mathematik und Informatik, Geb. 27, Universit\"at des Saarlandes,
D-66123 Saarbr\"ucken, Germany}\par
\noindent{schreyer@math.uni-sb.de}\par
}

\end{document}